\documentclass[a4paper,fleqn]{cas-sc}
\usepackage[authoryear,longnamesfirst]{natbib}
\usepackage[T1]{fontenc}
\usepackage[utf8]{inputenc}

\usepackage{amsmath,amssymb,amsfonts,amsthm}
\usepackage{booktabs}
\usepackage{array}
\usepackage{tabularx}
\usepackage{longtable}
\usepackage{multirow}
\usepackage{makecell}
\usepackage{adjustbox}
\usepackage[table]{xcolor}
\usepackage{colortbl}
\usepackage{graphicx}
\usepackage{float}
\usepackage{pdflscape}
\usepackage{etoolbox}
\usepackage{xspace}

\newtheorem{theorem}{Theorem}[section]
\newtheorem{lemma}[theorem]{Lemma}

\newcolumntype{P}{>{\centering\arraybackslash}p{0.7cm}}
\newcolumntype{Q}{>{\centering\arraybackslash}p{0.7cm}}

\begin{document}
	\let\WriteBookmarks\relax
	
	\shorttitle{Long-Lived Carpet-Like Transients}
	
	\shortauthors{M. Nowak-K\k{e}pczyk}
	
	% Title
	\title[mode=title]{Long-Lived Carpet-Like Transients in
		Time-Dependent Modular Discrete Laplacian Dynamics}

	% Author
	\author[1]{Ma\l gorzata Nowak-K\k{e}pczyk}[orcid=0000-0003-0934-9682]
	
	% E-mail
	\ead{[malnow@kul.pl](malnow@kul.pl)}
	%\ead{[gosianmk@wp.pl](gosianmk@wp.pl)}
	
	\credit{Conceptualization, Methodology, Software, Validation,
		Formal analysis, Investigation, Data curation, Visualization,
		Writing -- original draft, Writing -- review and editing}
	
	% Affiliation
	\affiliation[1]{
		organization={Institute of Informatics, John Paul II Catholic University of Lublin},
		addressline={Konstantynów 1H},
		city={Lublin},
		postcode={20-708},
		country={Poland}
	}
	%\cortext[1]{Corresponding author}
	\begin{abstract}		
		Binary modular Laplacian dynamics is organized by repeated replication, growth, and collapse on dyadic scales. We investigate how isolated and periodically repeated non-binary updates modify this organization and whether they can sustain densely occupied geometric structures.
		
		Computational experiments were performed for multiple finite seeds, neighborhood masks, inserted moduli, and periodic schedules. Constant prime moduli retain the expected \(p\)-adic replication hierarchy. A single non-binary insertion instead acts primarily as an effective-seed	replacement: subsequent evolution remains binary in character but may	resume at an integer phase shift of the original collapse--recovery	cycle, with only modest finite-time densification.
		
		A qualitatively different response occurs under periodic perturbation of binary regime, $[2,k,2^{s}]^\infty,$ where \(2^{s}\) denotes \(s\) consecutive binary updates.
		Odd-modulus insertions can sustain long-lived, spatially coherent carpet-like regimes with comparatively high and stable occupation, whereas even-modulus insertions remain  binary-like. Ternary insertion gives the broadest and most robust high-density response. Density depends	strongly and non-monotonically on the binary tail length; for \(k=3\), regularly spaced troughs near \(s=7,15,23\) reveal a phase-sensitive	interaction with the underlying binary epoch structure.
		
A striking cross-modulus synchronization also emerges. For the degree-four diagonal Neumann and von Neumann masks, the higher odd
moduli \(k=5,7,9\), although distinct during non-binary updates,
repeatedly collapse onto exactly the same complete configuration at
the following binary phases. A simple parity argument explains this loss of memory and the resulting coincidence.
		
		Finally, increasing seed extent generally raises occupation while attenuating tail-length sensitivity. The results identify periodic modular insertions as both a phase-selective and geometry-dependent mechanism capable of transforming \(p\)-adically organized replication
		into long-lived carpet-like dynamics.
		
	\end{abstract}
\maketitle


%\maketitle
%\justifying

%\maketitle

\section{Introduction}

Discrete Laplacians on infinite lattices form a natural class of
operators in discrete analysis, dynamical systems, and computational
mathematics. When iterated from finite seeds, they generate
cellular-automaton-like evolution rules that can be studied both
algebraically and through large-scale computational experiments.
Despite the simplicity of the underlying local interactions, the
resulting configurations may exhibit replication, symmetry
preservation, repeated collapse and rebuilding, and complex multiscale
geometry \cite{Ilachinski2001}.

A characteristic feature of modular discrete Laplacian dynamics is the
emergence of replication and fractal-like organization. Depending on
the modulus, seed, and neighborhood mask, the evolving configurations
may contain Sierpiński-type voids, replicated local templates, framed
structures, and densely perforated forms. Such behavior places the
subject naturally within experimental mathematics, where
computer-assisted exploration is used not only to visualize complex
dynamics, but also to formulate empirical classifications and
hypotheses concerning their temporal organization.

For constant prime moduli, part of this organization has a direct
algebraic origin. The Frobenius identity produces replication on
\(p\)-adic scales: finite templates reappear as translated copies,
while ordered sequences of growth, overlap, collapse, and rebuilding
may recur on progressively enlarged spatial and temporal scales. The
associated density and color entropy are therefore not generally
constant. Instead, they oscillate in relation to recurrent replication
and collapse stages.

Against this reference behaviour, we use the term
\emph{carpet-like regime} for a qualitatively different finite-time
organization. A carpet-like trajectory remains spatially coherent and
comparatively densely occupied, without the pronounced density
collapses and separation into sparse seed replicas characteristic of
constant-modulus dynamics. Its relative density remains comparatively
high and stable over the observation window, although the internal
texture may continue to evolve. The term is descriptive and refers
only to the investigated computational horizon.

The situation becomes less predictable when the modulus varies with
time. A single non-binary update may alter the active configuration,
produce a new effective seed, and displace the subsequent binary
trajectory without removing its basic collapse--recovery organization.
Periodically repeated interventions may have a stronger effect:
successive non-binary updates can interrupt the sparse replication
hierarchy and drive the system into the carpet-like regime defined
above.

\paragraph{Related work.}

Periodicity and replication phenomena in discrete Laplacian dynamics
were studied by \cite{Hadlich2011}. Fractal and geometric aspects of
related binary, ternary, quaternary, and higher-order discrete
structures were subsequently investigated in
\cite{Suzuki2018Geometrical,Lawrynowicz2019,NowakKepczyk2026}.

More broadly, the emergence of complex large-scale organization from
simple local interaction rules is a recurring theme in cellular
automata and self-organizing systems
\cite{Ilachinski2001,Whitesides2002,Singh2024NonEquilibrium}.
Although arising in a different physical setting, recent experimental
work has also demonstrated the spontaneous formation of fractal
geometries, including Sierpiński-like motifs, in protein assemblies
\cite{Sendker2024Fractal}. Such examples emphasize the broader
relevance of replication, symmetry inheritance, and multiscale
pattern formation generated by local rules.

\paragraph{Scope of the present study.}

We investigate modular iterations of discrete Laplacians in three
successive settings. First, constant-modulus dynamics is used to
establish the reference replication structure. Prime moduli are
interpreted through the Frobenius mechanism, while prime powers and
composite moduli are examined computationally to determine which
binary-like, ternary-like, or mixed signatures persist beyond the
direct finite-field argument.

This baseline also makes it possible to distinguish the roles of the
initial seed, neighborhood mask, and modular schedule. The mask
primarily determines the large-scale envelope of the evolving
configuration, while shared seed--mask symmetries constrain its
preserved reflection symmetries. The modulus and its temporal schedule
control the internal morphology, density evolution, replication clock,
and sequence of collapse and rebuilding stages.

Second, we examine isolated non-binary insertions into an otherwise
binary evolution. We ask whether a single intervention destroys the
binary organization or instead replaces the effective seed and shifts
the subsequent trajectory. This provides a bridge between the
constant-modulus baseline and genuinely periodic perturbations.

Finally, we study the periodic family
\begin{equation}
	\sigma_{k,s}
	=
	[2,k,2^{s}]^\infty,
	\label{eq:intro-periodic-family}
\end{equation}
where \(2^{s}\) denotes \(s\) consecutive updates modulo \(2\).
The inserted modulus \(k\) determines the type of intervention, whereas
the tail length \(s\) controls the time available for binary
development before the next intervention. The shortest schedules,
\[
[2,k]^\infty,\qquad
[2,k,2]^\infty,\qquad
[2,k,2,2]^\infty,
\]
provide cross-sections of the transition from maximal alternating
perturbations to progressively longer binary intervals.

For the complete family, we quantify mean density and temporal
variability over common finite observation windows. Particular
attention is given to the contrast between even- and odd-modulus
interventions, the non-monotone dependence on binary tail length, and
the unusually robust response produced by ternary insertions. We also
examine how increasing seed extent and using a more overlap-promoting
mask can attenuate the visible sensitivity to the modular schedule.

\section{Model and definitions}

We consider time-dependent modular dynamics on the infinite square
lattice $\mathbb Z^2$. The evolution is generated by successive
applications of a discrete Laplacian operator, while the modulus used
for the reduction may vary with time.

We treat a non-binary update introduced into the binary dynamics as a
local perturbation of the modular schedule. Throughout the paper, such
an update is referred to as a \emph{non-binary insertion}, or simply an
\emph{insertion}.

\subsection{Modular Laplacian dynamics}

A configuration at time \(t\) is a function
\[u_t:\mathbb Z^2\longrightarrow\mathbb Z_{\geq 0}\]
with finite support. Its support is
\[A_t=\operatorname{supp}(u_t)
=\{p\in\mathbb Z^2:u_t(p)\neq0\}.\]

Let \(M\subset\mathbb Z^2\setminus\{(0,0)\}\) be a fixed finite
neighborhood mask. For a lattice site \(p\in\mathbb Z^2\), let
\[
N_M(p)=\{p+v:v\in M\},
\qquad
d_M=|M|.
\]

The mask Laplacian is obtained by summing, over all active positions
of the mask, the difference between the neighboring value and the
value at the central site:
\begin{equation}
	(L_Mu)(p)=\sum_{q\in N_M(p)}
	\bigl(u(q)-u(p)\bigr)
	=
	\sum_{q\in N_M(p)}u(q)-d_Mu(p).
	\label{integer_laplacian}
\end{equation}

For an integer \(x\) and a modulus \(k\geq2\), let
\([x]_k\in\{0,1,\ldots,k-1\}\) denote the standard representative of
\(x\) modulo \(k\). The evolution is defined by
\begin{equation}
	u_t(p)=
	\bigl[(L_Mu_{t-1})(p)\bigr]_{k_t},
	\qquad t\geq1,
	\label{modular_evolution}
\end{equation}
where $u_0$ is the initial configuration and
\[k_1,k_2,k_3,\ldots,\qquad k_t\in\{2,3,4,\ldots\},
\label{sequence}\]
is the prescribed sequence of moduli. Thus, the first update from
$u_0$ to $u_1$ is performed modulo $k_1$.
The constant schedule $(2,2,2,\ldots)$
will be referred to as the binary dynamics. Analogously, constant
updates modulo $3$, $4$, and higher moduli will be referred to as
ternary, quaternary, and higher-order dynamics, respectively.

\subsection{Modular schedules}

Let
\[
B=(b_1,b_2,\ldots,b_P)
\]
be a finite block of moduli. Its periodic repetition is denoted by
\[
[b_1,b_2,\ldots,b_P]^\infty,
\]
meaning that
\[
k_t=b_{\,1+((t-1)\bmod P)}, \qquad t\geq 1.
\]
Thus, square brackets followed by the superscript \(\infty\) always
denote periodic repetition of the displayed block. In particular,
\([k]^\infty\) denotes constant-modulus dynamics.

A single non-binary insertion into an otherwise binary schedule is
written as
\[\eta_k=(2,k,2,2,2,\ldots), \qquad k\neq 2.\]
Hence \(k_1=2\), \(k_2=k\), and \(k_t=2\) for all \(t\geq3\).
This schedule is not periodic.

Periodic insertions are described by the family
\[\sigma_{k,s}=[2,k,2^{s}]^\infty,\qquad s\geq0,\]
where \(2^{s}\) denotes \(s\) consecutive copies of the modulus
\(2\), not the numerical power \(2^s\). Thus,
\[
\sigma_{k,0}=[2,k]^\infty,\qquad
\sigma_{k,1}=[2,k,2]^\infty,\qquad
\sigma_{k,2}=[2,k,2,2]^\infty.
\]

The period of \(\sigma_{k,s}\) is
\[P=s+2,\]
and successive updates modulo \(k\) are separated by \(s+1\) binary
updates. We refer to \(k\) as the \emph{inserted modulus} and to \(s\)
as the \emph{binary-tail length}. The specific values of \(k\) and
\(s\) used in each experiment are given in Section~2.5.

\subsection{Seeds and neighborhood masks}

The initial configuration $u_0$ is called a \emph{seed}. Throughout
the study, seeds have finite support.

A \emph{figure} associated with a configuration $u_t$ is its geometric
support \[F_t=A_t.\]

In earlier studies of discrete Laplacian dynamics, seeds were typically
contained in a $3\times3$ square. We refer to such configurations as
\emph{small seeds}. Seeds contained in squares of size $18\times18$
and $84\times84$ are referred to as \emph{medium} and \emph{large seeds}, respectively. These experimental size classes represent localized, intermediate, and large-scale initial configurations while remaining computationally tractable over the investigated time
horizons.

Every finite figure may itself be used as a 
seed. Small, medium, and large seed classes are used as reference
categories in the computational experiments.

%Fig=1
\begin{figure}[pos=htbp]%Fig1
	\centering
	\includegraphics[width=0.8\textwidth]{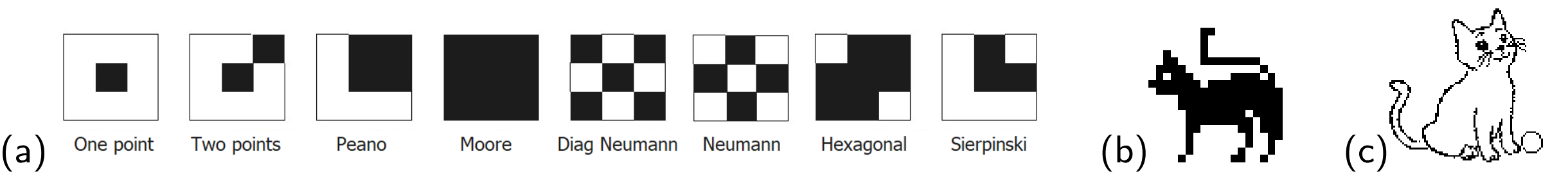}
	\caption{Examples of experimental seed classes:
		(a) small, (b) medium, and (c) large seeds.}
	\label{seed}
\end{figure}

\paragraph{Neighborhood masks.}

Each neighborhood mask is represented by a finite set
\[M\subset\mathbb Z^2\setminus\{(0,0)\}.\]
The same mask $M$ is used at every iteration of a given evolution.
Only the modulus $k_t$ may change with time.

\begin{figure}[pos=htbp]%Fig2
	\centering
	\includegraphics[width=0.64\textwidth]{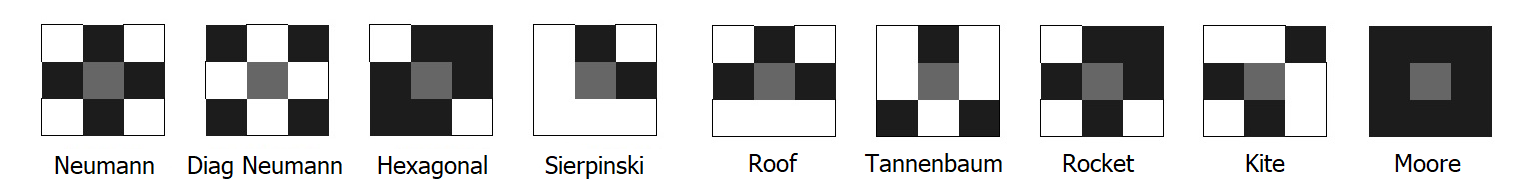}
	\caption{Examples of neighborhood masks used in the computational
		experiments. The central cell is not included in the neighborhood.}
	\label{neighbor}
\end{figure}

We consider several neighborhood masks, including the von Neumann,
diagonal Neumann, Moore, Tannenbaum, hexagonal, and other
neighborhoods shown in Fig.~\ref{neighbor}. Less common masks are
included to examine how local neighborhood geometry affects the
large-scale envelope and symmetry of the resulting configurations.

%Fig=2

%%%%%%%%%%%%%%%%%%%%%%%%%%%%%%%%%%%%%%%%%%%%%%%%%%%
\subsection{Observables}

\paragraph{Replication and recurrence.}

Replication and revival phenomena for constant-modulus Laplacian
dynamics, including overlapping, spatially separated, and shifted
replicas, were characterized in detail in our previous work
\cite{NowakKepczykFrobeniusRevivals}. In the present study, these
phenomena are used only as reference behavior for the binary and
constant-modulus dynamics and are not reclassified.

Recurrence was tested at the level of complete configurations,
including their residue values. An exact recurrence was recorded when
there existed \(t_0\geq 0\) and \(\tau>0\) such that
\[u_{t_0+\tau}(p)=u_{t_0}(p)\qquad\text{for every }p\in\mathbb Z^2.\]
A translated recurrence was recorded when there additionally existed
\(v\in\mathbb Z^2\) such that 
\[u_{t_0+\tau}(p)=u_{t_0}(p-v)
\qquad\text{for every }p\in\mathbb Z^2.\]
Unless stated otherwise, the absence of recurrence means that neither
condition was detected within the investigated simulation horizon.

\paragraph{Relative density and finite-window statistics.}

Let
\[A_t=\{p\in\mathbb Z^2:u_t(p)\neq0\}\]
be the support of the configuration at iteration \(t\), and let
\[P_t=\operatorname{conv}(A_t)\]
be the smallest convex polygon containing \(A_t\). We refer to \(P_t\)
as the bounding polygon of the configuration.

Let $
\Lambda_t=P_t\cap\mathbb Z^2$
denote the set of lattice sites contained in this polygon. For a
nonempty configuration, the relative density is defined as
\begin{equation}
	\rho_t=
	\frac{|A_t|}{|\Lambda_t|}.
	\label{relative_density}
\end{equation}

If \(A_t=\varnothing\), we set \(\rho_t=0\).

This definition evaluates each configuration relative to its natural
polygonal envelope, independently of its orientation.

Let
$W=\{t_a,t_a+1,\ldots,t_b\}$
be a finite evaluation window. The mean density in \(W\) is
\begin{equation}
	\mu_W=\frac{1}{|W|}	\sum_{t\in W}\rho_t.
	\label{mean_density_window}
\end{equation}

The density standard deviation is
\begin{equation}
	SD_W=\sqrt{\frac{1}{|W|-1}\sum_{t\in W}	\left(\rho_t-\mu_W\right)^2	}.
	\label{sd_density_window}
\end{equation}

Relative density variability is quantified by the coefficient of
variation
\begin{equation}
	CV_W=\frac{SD_W}{\mu_W},\qquad\mu_W>0.	\label{cv_density_window}
\end{equation}

Mean density and the coefficient of variation are treated as separate
continuous characteristics. High \(\mu_W\) indicates that a relatively
large proportion of the bounding polygon remains occupied, whereas
low \(CV_W\) indicates that this density is maintained comparatively
regularly throughout the evaluation window.

\paragraph{Color entropy.}

For a configuration obtained at time \(t\) using modulus \(k_t\), let
\[
\pi_t^{(c)}=
\frac{\left|\{p\in \Lambda_t : u_t(p)=c\}\right|}{|\Lambda_t|},
\qquad c=0,1,\ldots,k_t-1,
\]
where \(\Lambda_t\) denotes the lattice sites contained in the geometric
envelope used to calculate relative density.

The color entropy is defined as
\begin{equation}
	H_t=-\sum_{c=0}^{k_t-1}
	\pi_t^{(c)}\log \pi_t^{(c)},
	\qquad 0\log 0:=0.
	\label{color_entropy}
\end{equation}
For diagnostic purposes, the contribution of residue class \(c\) to
the total colour entropy is denoted by
\[
h_t^{(c)}=-\pi_t^{(c)}\log \pi_t^{(c)},
\]
so that \(H_t=\sum_{c=0}^{k_t-1} h_t^{(c)}\).

Color entropy describes the distribution of residue values, whereas
relative density distinguishes only zero and nonzero cells. It is used
as a supplementary measure, particularly when comparing dynamics under
different moduli and when examining composite-modulus cases such as
\([6]^\infty\).

\paragraph{Exact configuration coincidence.}
Configurations generated for different inserted moduli were classified
as identical only when
\[
u_t^{(k_1)}(x)=u_t^{(k_2)}(x)
\qquad\text{for all }x\in\mathbb Z^2.
\]
Comparisons were made at matched global iterations for the same seed,
mask, and short-tail schedule.

\subsection{Experimental design}

\paragraph{Exploratory phase.}
Initial computations were performed across a broader range of finite
seeds, neighborhood masks, moduli, and modular schedules. These runs
were used to identify recurrent geometric and dynamical effects,
assess their dependence on seed and mask geometry, and formulate the
specific questions examined below. They were exploratory rather than
a balanced factorial experiment; quantitative results reported in the
paper refer to the parameter settings specified for the corresponding
experiments.

\paragraph{Constant-modulus experiments.}
The constant schedules
\[
[k]^\infty=(k,k,k,\ldots),
\qquad k\in\{2,3,4,5,6,7,8,9\},
\]
were examined as reference dynamics. Several seeds and neighborhood
masks were used during the exploratory comparison to identify
replication, shape-inheritance, and scalar-oscillation signatures.
Representative cases from these computations are reported in
Section~3.2 for prime, prime-power, and composite moduli.

\paragraph{Single-insertion experiments.}
We evaluated the schedules \(\eta_k\) defined in Section~2.2 for
\[
k\in\{3,4,5,6,7,9\}.
\]
The experiments examined whether a single non-binary insertion changed
the effective seed or the phase of the subsequent binary dynamics.
Density effects were summarized over \(t=2,\ldots,82\), and phase
comparisons with the binary reference were performed over
\(t=17,\ldots,80\) for the reported seed--mask combinations.

\paragraph{Periodic-insertion experiments.}
We evaluated the schedules \(\sigma_{k,s}\) defined in Section~2.2 for
\[
k\in\{3,4,5,6,7,9\}.
\]
Short-tail schedules with \(s=0,1,2\) were first used to examine the
transition between alternating insertion and increasing binary
recovery. The resulting modulus classes motivated a broader
investigation of binary-tail length. Additional computations varied
seed extent and neighborhood geometry to assess the robustness and
geometric dependence of the observed effects. Representative results
are reported in Section~3.4; the parameter settings for individual
analyses are specified with the corresponding results and figures.

\paragraph{Simulation and evaluation.}
Unless stated otherwise, periodic trajectories were simulated up to
\[T=200.\]
The computational domain was expanded at every iteration by the radius
of the neighborhood mask, preventing truncation at a fixed numerical
boundary. The first \(32\) iterations were treated as a developmental
stage, and periodic-schedule statistics were evaluated over
\[W=[33,200].\]

\paragraph{Outcome measures.}
The principal outcome measures were the mean relative density
\(\mu_W\) and the coefficient of variation of density \(CV_W\).
These quantities describe complementary properties of the trajectory:
\(\mu_W\) measures its average spatial occupation, whereas \(CV_W\)
measures relative temporal variability.

\paragraph{Aggregation and comparison.}

For the \(s=2\) comparison, density statistics were summarized across
the seed--mask panel using medians and interquartile ranges. The
systematic tail-length scan was then performed for three small seeds
with the diagonal Neumann mask, the setting in which schedule effects
were most clearly resolved. For each \((k,s)\), median mean density
across the three seeds and the corresponding inter-seed IQR were
reported. Seed-size effects were examined separately.

\paragraph{Reproducibility.}
The source code used for the computational experiments and the input
seed configurations are available in the accompanying Figshare
repository. The individual scripts specify the seed, neighborhood mask, modular schedule, simulation horizon, and evaluation settings used in the corresponding analyses. 

Generative AI assistance (ChatGPT, OpenAI) was used during the
review and debugging of selected Python scripts; all code and
computational outputs were subsequently inspected, tested, and
validated by the author.

\subsection{Algebraic background}\label{subsec:algebraic-background}

The replication phenomena observed under constant prime moduli are
related to the Frobenius endomorphism and to the vanishing of mixed
multinomial terms in finite characteristic.

Let \(\tau_v\) denote the lattice translation operator
\[(\tau_v u)(p)=u(p+v),\qquad v\in\mathbb Z^2.\]
For a fixed neighborhood mask \(M\), the discrete Laplacian can be
written as
\begin{equation}
	L_M	=\sum_{v\in M}\tau_v-d_M I,\qquad d_M=|M|.\label{translation_laplacian}
\end{equation}

All translation operators commute. Let \(p\) be prime and consider
the dynamics over the finite field \(\mathbb F_p\). In characteristic
\(p\), the Frobenius identity gives
\[(a+b)^{p^m}
=a^{p^m}+b^{p^m}.\]
Its multinomial form therefore yields
\begin{equation}
	\left(\sum_j a_j T_j\right)^{p^m}
	=\sum_j a_j^{p^m}T_j^{p^m}
\end{equation}
for mutually commuting operators \(T_j\).

Applying this identity to \(L_M\), and using
$\tau_v^{\,p^m}=\tau_{p^m v},$
gives
\begin{equation}
	L_M^{p^m}=\sum_{v\in M}\tau_{p^m v}	-
	d_M I\pmod p,\label{frobenius_laplacian}
\end{equation}
where the coefficients are interpreted in \(\mathbb F_p\).

Thus, at the characteristic scales \(p^m\), the mixed operator terms
vanish and the evolution is assembled from spatially rescaled
translation components. When these translated contributions do not
overlap or cancel, the support contains replicated copies of the
earlier configuration.

The characteristic algebraic scales are therefore
\[2^m=1,2,4,8,\ldots\]
for binary dynamics,
\[3^m=1,3,9,27,\ldots\]
for ternary dynamics, and
\[5^m=1,5,25,125,\ldots\]
for dynamics modulo \(5\), with analogous relations for other prime
moduli.

These characteristic scales should be distinguished from the minimal
support-replication times observed for a particular seed and mask.
The latter additionally depend on the geometry of the initial
configuration, possible overlap or cancellation of translated
components, and the separation condition for their bounding
polygons.

For composite moduli and prime powers, the finite-field Frobenius
identity does not apply directly to the full dynamics over
\(\mathbb Z/k\mathbb Z\). Binary-like behavior observed for moduli
\(4\) and \(8\), ternary-like behavior observed for modulus \(9\), and
mixed behavior for other composite moduli are therefore treated below
as computational observations rather than direct consequences of
Eq.~\eqref{frobenius_laplacian}.

\section{Results}
%%%%%%%%%%%%%%%%%%%%%%%%%%%%%%%%%%%%%%%%%%%%%%%%%%%%%%%%%
%\FloatBarrier
\subsection{Mask-induced geometry and symmetry}
\label{subsec:shape-inheritance}

\begin{figure}[pos=htbp]
	\centering
	\includegraphics[width=0.9\textwidth]{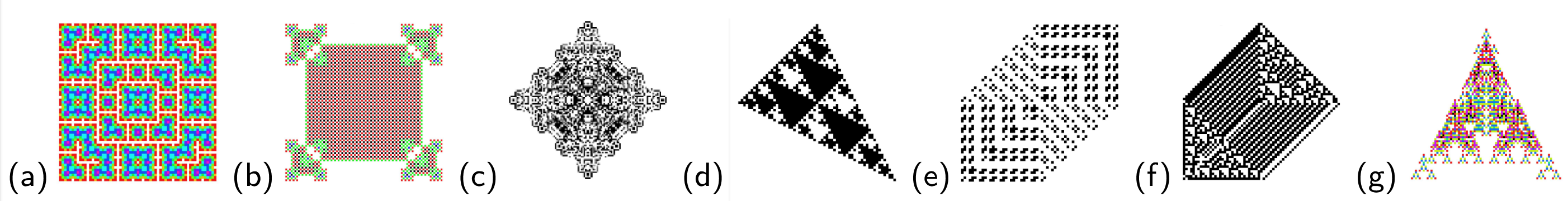}
	\caption{
		\textbf{Shape inheritance across neighborhood masks.}
		(a) Moore and (b) diagonal Neumann masks produce square-type
		envelopes; (c) von Neumann produces a diamond; (d) kite, a skew
		triangle; (e) hexagonal, a hexagon; (f) rocket, a pentagon-like
		envelope; and (g) Tannenbaum, a triangle.
		All examples were generated from the same medium seed under
		constant-modulus dynamics.
	}
	\label{fig:inheritance}
\end{figure}

Across the tested seeds, moduli, and schedules, the large-scale
envelope reflected the geometry of the neighborhood mask. We refer to
this property as \emph{shape inheritance}. Moore and diagonal Neumann
masks produced square-type envelopes, the von Neumann mask produced a
diamond-shaped envelope, and the remaining masks generated skew
triangular, hexagonal, pentagon-like, or triangular outlines
(Fig.~\ref{fig:inheritance}).

%Fig=3

The outer envelope is determined primarily by the mask, whereas global
reflection symmetry depends on both the seed and the mask. Let \(G_S\)
and \(G_M\) denote their symmetry groups. Since the Laplacian update is
equivariant under every symmetry shared by the two, the symmetries in
\[
G_S\cap G_M
\]
are preserved throughout the evolution. No additional global
reflection symmetry is guaranteed.

The modular schedule mainly affects the internal morphology and its
temporal development, including growth, collapse, recovery, and the
formation of repeated holes or framed structures. In some cases, these
holes form a negative Sierpiński-type organization: recursively arranged
voids appear within an expanding occupied region, as in
Fig.~\ref{fig:inheritance}(d).
	
%%%%%%%%%%%%%%%%%%%%%%%%%%%%%%%%%%%%%%%%%%%%%%%%%%%%%%%%%%%%%%%%

\subsection{Constant-modulus dynamics}
\label{subsec:constant-modulus}

Constant schedules provide the baseline for the time-dependent dynamics
studied below:
\[
[k]^\infty=(k,k,k,\ldots).
\]

\subsubsection{Prime moduli: Frobenius scales and \(p\)-adic recurrence}

For a prime modulus \(p\), the Frobenius identity selects the
characteristic scales
\[
p^m,\qquad m=0,1,2,\ldots .
\]
At these scales, mixed operator terms vanish and the configuration is
reorganized into translated copies of an earlier state.

This replication mechanism has two geometric consequences. First,
finite templates, including the seed itself, return as translated
copies at \(p\)-adic separations. Second, the repeated separation,
growth, and overlap of these copies generates recurrent oscillations
of support density and color entropy. The exact oscillation amplitudes
depend on the seed and mask, whereas their characteristic organization
is tied to the \(p\)-adic replication clock.

\paragraph{Local \(p\)-adic recurrence.}
For a fixed local state \(u_r\), the Frobenius identity gives
\[
u_{p^m+r}
=
\sum_{\nu}a_\nu T_{p^m\nu}u_r.
\]
Thus, whenever the translated components remain separated, the same
finite template reappears as spatially translated copies at
\(p\)-adic separations. The copies may differ by multiplication by a
non-zero residue, but their support and local morphology are
preserved. We refer to this property as \emph{local \(p\)-adic
	recurrence}; for \(p=2\), it becomes local dyadic recurrence.

%Fig=4
\begin{figure}[pos=htbp]
	\centering
	\includegraphics[width=0.8\textwidth]{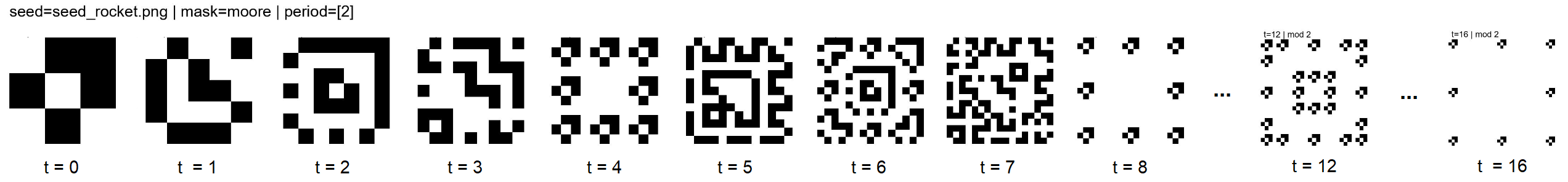}
	\caption{The panels show successive seed-return stages. Although their displayed
		iteration times need not themselves be powers of two, the number and
		mutual separation of the returned templates are organized by the
		underlying dyadic replication hierarchy.
	}
	\label{fig:local_dyadic_seed_return}
\end{figure}

\begin{figure}[pos=htbp]
	\centering
	\includegraphics[width=\textwidth]
	%{binary_epochs_three_seeds.png}
	{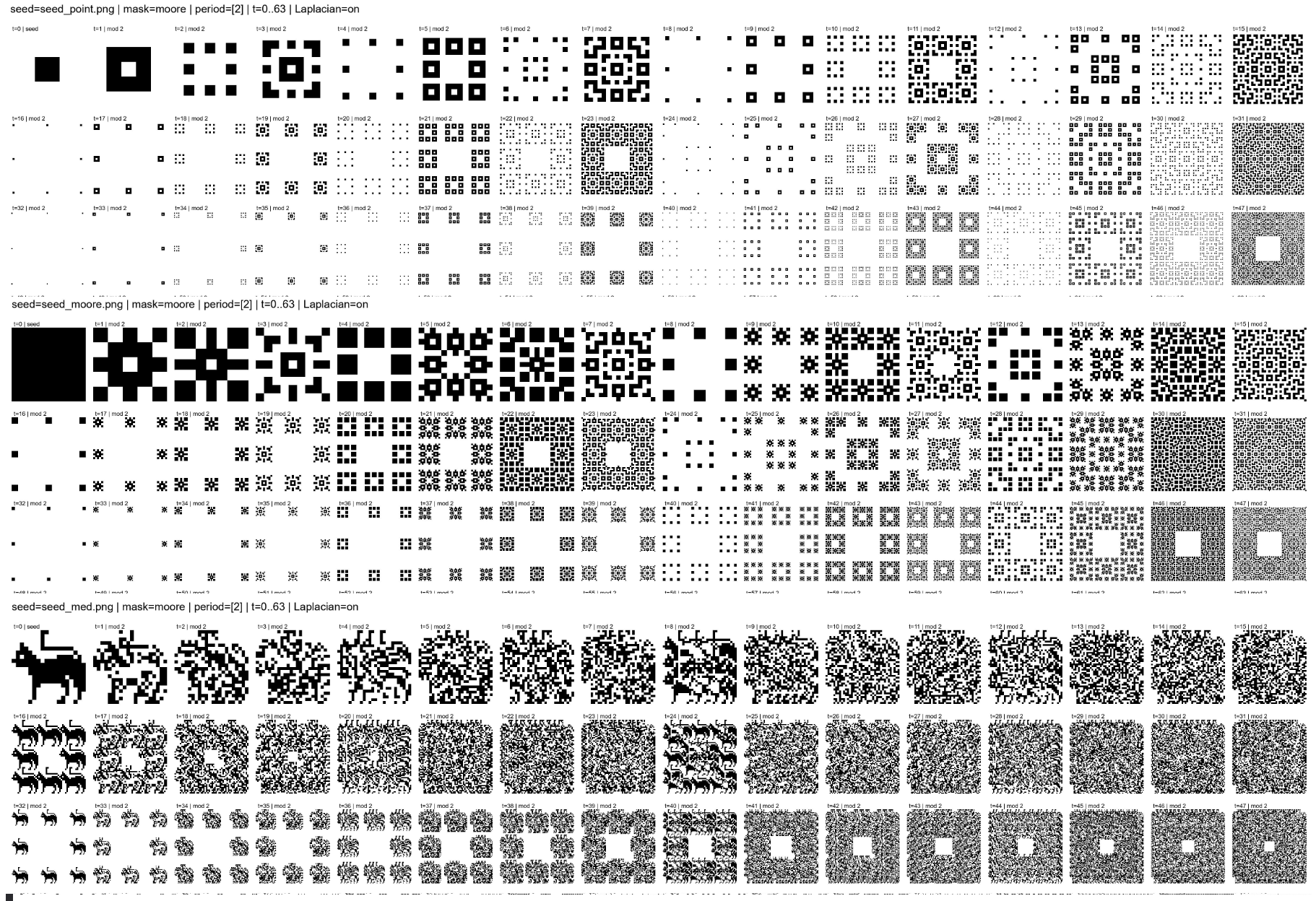}
	\caption{\textbf{Global dyadic recurrence in binary dynamics.}
		Point, Moore, and medium image seeds evolved with the Moore mask for
		\(t=0,\ldots,63\). All three pass through corresponding two-part
		shape epochs with lengths \(4,4\), \(8,8\), and \(16,16\), followed
		by collapse into increasingly separated seed copies. Larger seeds
		produce longer overlap but preserve the same epoch structure.}
	\label{fig:binary_seed_robust_epochs}
\end{figure}

\begin{figure}[pos=htbp]
	\centering
	\includegraphics[width=0.8\textwidth]
	%{ternary_epoch_recurrence.png}
	{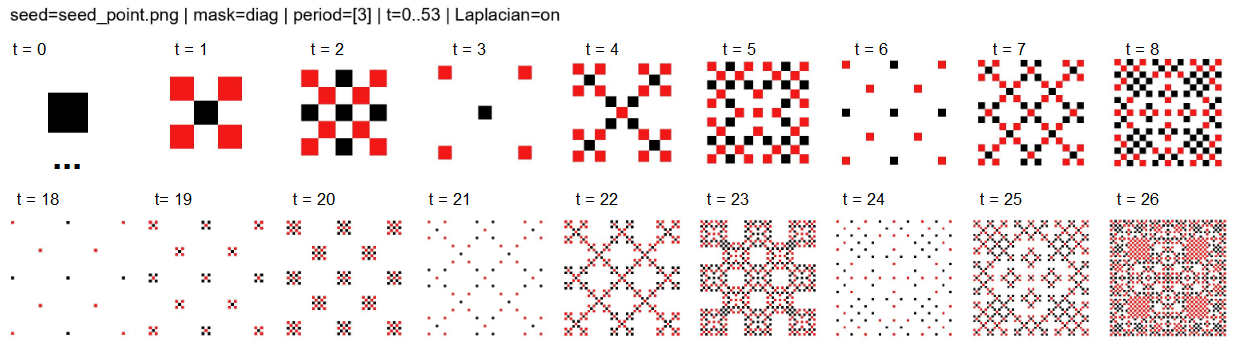}
	\caption{\textbf{Local and global ternary recurrence under
			constant modulo-\(3\) dynamics.}
		Point-seed evolution with the diagonal Neumann mask under
		\([3]^\infty\). The rows show the reference sequence
		\(u_0,\ldots,u_8\) and its later realizations
		\(u_{18},\ldots,u_{26}\).
	Within the later sequence, three recurring three-step motif families
	are visible: a sparse dot array, a skew lattice, and a compact
	four-petal pattern.}
	\label{fig:ternary_epoch_recurrence}
\end{figure}

\paragraph{Global \(p\)-adic epoch repetition.}
The same replication operator acts not only on an isolated seed-return
state, but on an ordered sequence of states \(u_r\). Consequently, an
entire seed--growth--overlap--collapse sequence may reappear in a
larger configuration as translated copies organized on the scale
\(p^m\). We call this recurrence of the ordered developmental sequence
\emph{global \(p\)-adic epoch repetition} (see Figs.~\ref{fig:binary_seed_robust_epochs}, \ref{fig:ternary_epoch_recurrence}). The global configuration
continues to expand and therefore does not return in fixed spatial
coordinates.

\paragraph{Binary dynamics.}
Binary evolution alternates between growth, collapse, and dispersal into
separated seed copies. Between the main collapse events it passes
through repeated two-part shape epochs. In
Fig.~\ref{fig:binary_seed_robust_epochs}, the corresponding phase
lengths are
\[
4,4,\;8,8,\;16,16,
\]
so the same two-stage motif sequence returns on a doubled scale.

The point, Moore, and medium image seeds all follow this organization.
They differ in local texture and in the duration of overlap, but not in
the main sequence of motif formation, collapse, and replication. Larger
seeds simply hide the scaffold for longer. We therefore use the point
seed as the clearest probe of the operator-driven dynamics, while
retaining larger seeds for selected robustness checks.

%Fig=5

\paragraph{Further prime moduli.}
The ternary case in Fig.~\ref{fig:ternary_epoch_recurrence} shows that
the distinction between local recurrence of finite templates and
global recurrence of an ordered developmental sequence is not specific
to the binary rule. The same Frobenius mechanism applies to further
prime moduli: the dyadic clock \(2^m\) is replaced by the corresponding
\(p\)-adic clock. Computations for \(p=5\) and \(p=7\) likewise showed
recurrent growth, collapse into separated seed copies, and rebuilding
on larger scales.

%Fig=6

\begin{figure}[pos=htbp]
	\centering
	\includegraphics[width=0.7\textwidth]{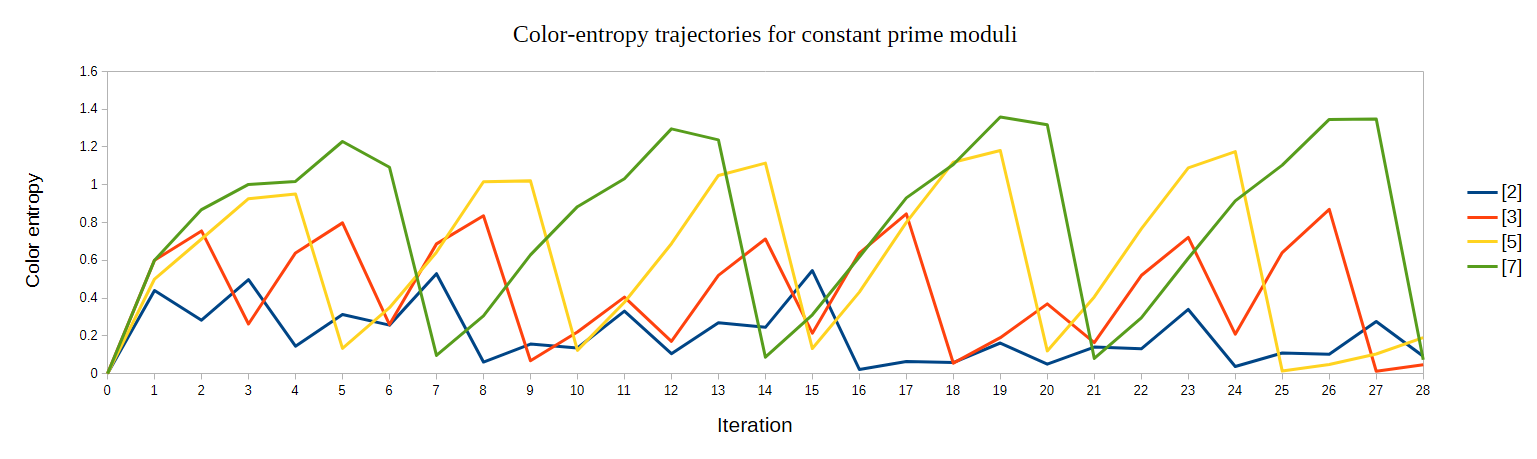}
	\caption{\textbf{Color-entropy trajectories for constant
			prime-modulus dynamics.}
		Color entropy \(H_t\) for the point seed and diagonal Neumann mask
		under \([2]^\infty\), \([3]^\infty\), \([5]^\infty\), and
		\([7]^\infty\). Recurrent entropy decreases are synchronized with
		sparse replication and collapse stages. Since the maximum raw entropy
		equals \(\log p\), the figure is intended to display temporal
		organization rather than absolute entropy levels.}
	\label{fig:Color_entropy_primes}
\end{figure}

\paragraph{Scalar signatures of \(p\)-adic recurrence.}
The geometric recurrence is accompanied by oscillations of scalar
observables. Support density increases during growth and overlap and
decreases when the configuration separates into sparse replicas.
Color entropy follows a related pattern: residue diversity rises during
motif development and falls at sparse replication stages
(Fig.~\ref{fig:Color_entropy_primes}). These observables do not remain
constant, but their oscillations are organized by the characteristic
\(p\)-adic recurrence scales.

%Fig=7

\subsubsection{Prime-power moduli}

Moduli such as \(4\), \(8\), and \(9\) are not finite fields, so the
Frobenius argument does not apply directly to the full dynamics over
\(\mathbb{Z}/k\mathbb{Z}\). Nevertheless, the computations showed
clear analogies:
\begin{itemize}
	\item \(4\) and \(8\): binary-like signatures,
	\item \(9\): ternary-like signatures.
\end{itemize}
These similarities are reported as computational observations rather
than direct algebraic consequences.

\subsubsection{Composite moduli: mixed prime-factor signatures}

The composite modulus \(6=2\cdot3\) exhibits a mixed temporal
organization. This structure is particularly clear in the colour
entropy, where fluctuations occur at characteristic scales associated
with both the binary and ternary components. The observation is
consistent with the Chinese-remainder decomposition
\(\mathbb Z_6\cong\mathbb Z_2\times\mathbb Z_3\): the dynamics retain
signatures of both constituent prime-modulus systems rather than
forming a simple intermediate regime.

%Fig=8
\begin{figure}[pos=htbp]%Fig13
	\centering
	\includegraphics[width=0.8\textwidth]{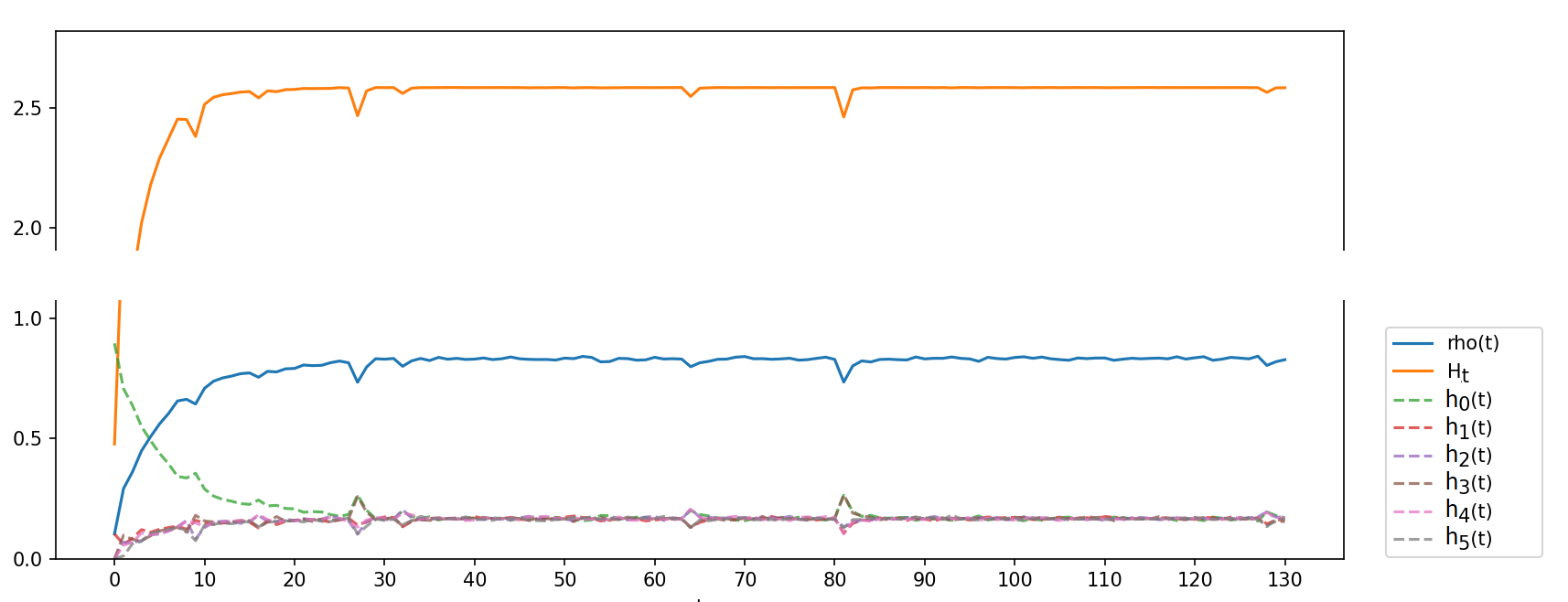}
\caption{
	Colour-entropy dynamics for the constant-modulus sequence
	\([6]^\infty\). Entropy fluctuations occur at temporal scales
	characteristic of both the binary and ternary components, revealing
	the mixed organization expected from
	\(\mathbb Z_6\cong\mathbb Z_2\times\mathbb Z_3\).
}
	\label{fig:Hc-6}
\end{figure}

\begin{figure}[pos=htbp]
	\centering
	\includegraphics[width=0.25\textwidth]{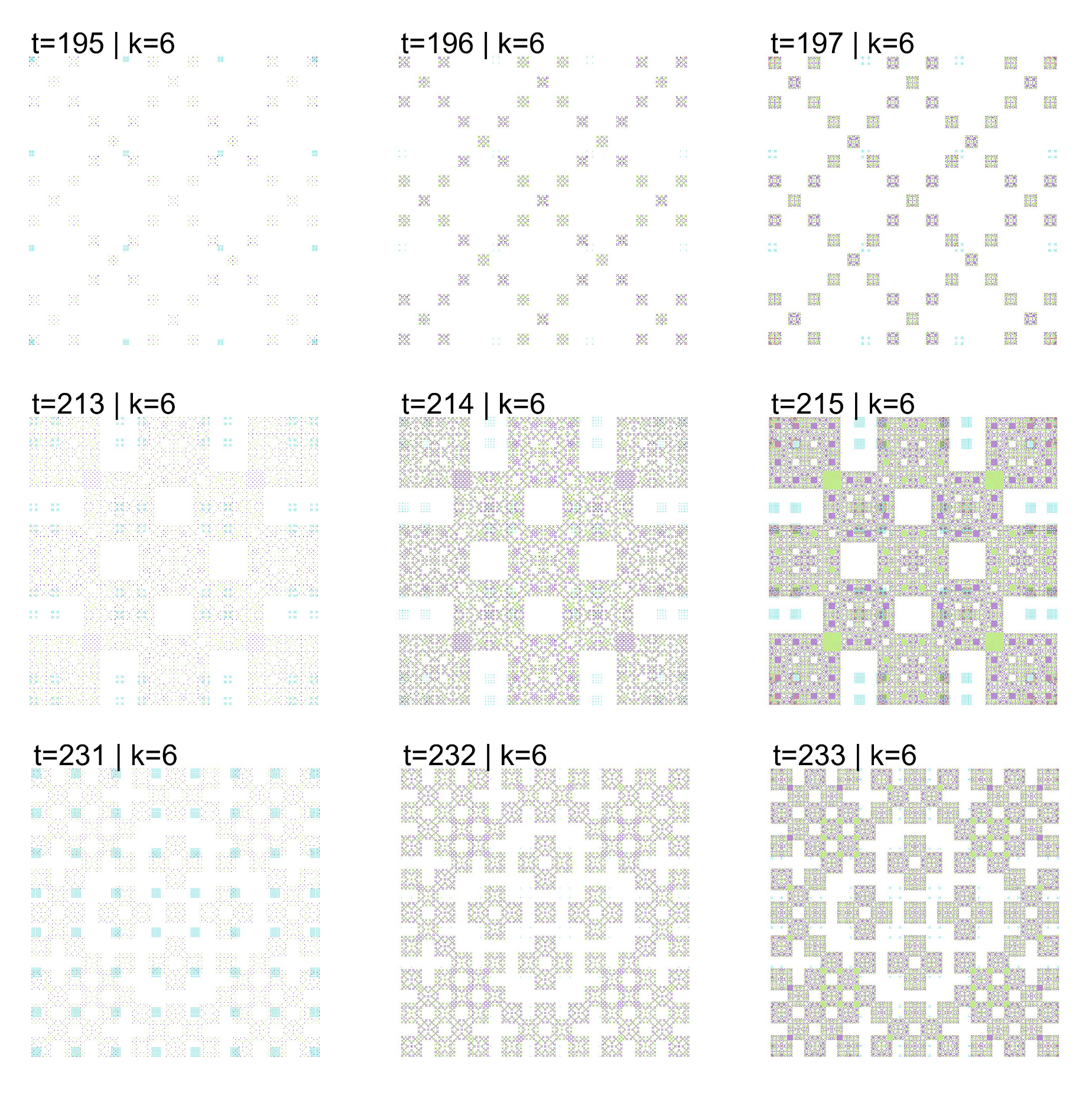}
	$\;\;\;\;$
	\includegraphics[width=0.25\textwidth]{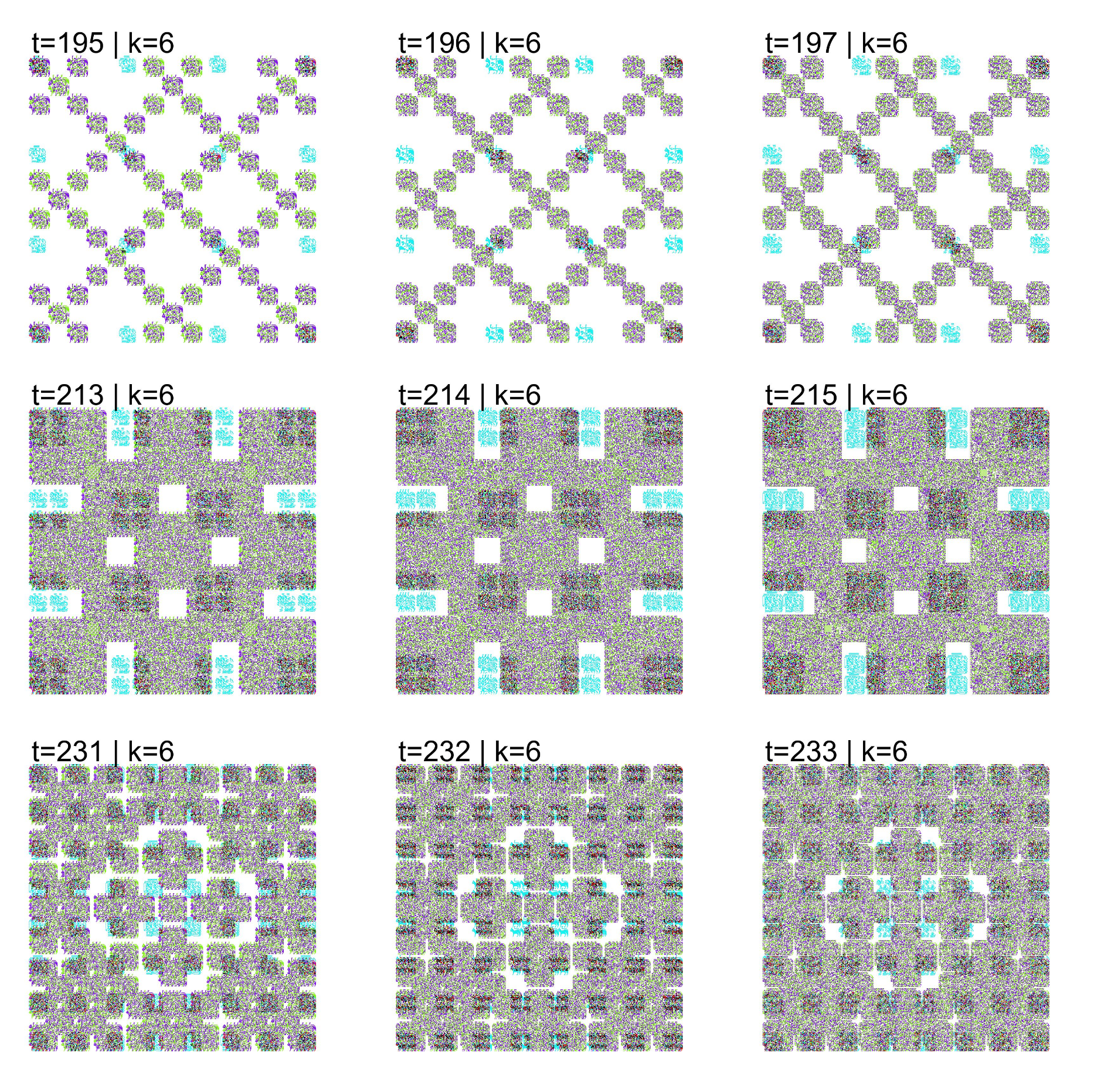}
	\caption{\textbf{Corresponding mesoscopic motifs for different seeds
			under constant modulo-\(6\) dynamics.}
		Point seed (left) and medium cat seed (right), evolved with the
		diagonal Neumann mask under \([6]^\infty\). The same mesoscopic
		motif sequence is visible in both cases. The image seed fills the
		operator-generated scaffold with a larger local glyph, prolonging
		overlap between neighbouring components.}
	\label{fig:mod6_seed_comparison}
\end{figure}

Despite this mixed temporal structure, corresponding mesoscopic motifs
remain visible for different seeds. Figure~\ref{fig:mod6_seed_comparison}
compares the point seed with the medium cat seed. The same diagonal
lattice-like, framed-square, and compact motifs appear at corresponding
iterations. The point seed reveals the scaffold directly, whereas the
image seed populates it with translated copies of the original motif.

%Fig=9

\paragraph{Summary of constant moduli}

The principal constant-modulus observations are summarized in
Table~\ref{tab:constant_moduli_summary}.

\begin{table}[pos=htbp]
	\centering
	\caption{
		Summary of recurrence and scalar-oscillation signatures observed
		for representative constant-modulus dynamics. For prime moduli,
		the \(p\)-adic replication scales follow from the finite-field
		argument, whereas the detailed geometry and the associated density
		and entropy oscillations are characterized computationally.
	}
	\label{tab:constant_moduli_summary}
	
	\begin{tabularx}{\textwidth}{
			>{\raggedright\arraybackslash}p{1.5cm}
			>{\raggedright\arraybackslash}p{7.2cm}
			X
		}
		\toprule
		\textbf{Schedule}
		& \textbf{Characteristic dynamical signature}
		& \textbf{Basis and status}
		\\
		\midrule
		
		\([2]^\infty\)
		&
		Local dyadic recurrence of seed templates and global dyadic
		epoch repetition of the seed--growth--overlap--collapse cycle;
		oscillatory density and entropy organized on scales \(2^m\)
		&
		Finite-field replication mechanism; recurrence geometry and
		scalar oscillations documented computationally
		\\
		
		\([3]^\infty\)
		&
		Local and global ternary recurrence; recurrent collapse and
		rebuilding accompanied by oscillatory density and entropy on
		scales \(3^m\)
		&
		Finite-field replication mechanism; detailed recurrence
		illustrated computationally
		\\
		
		\([5]^\infty,\,[7]^\infty\)
		&
		\(p\)-adic replication, recurrent seed-return and rebuilding,
		with density and entropy oscillations organized by the
		corresponding scales \(p^m\)
		&
		Finite-field replication mechanism; geometric and scalar
		signatures observed computationally
		\\
		
		\([4]^\infty,\,[8]^\infty\)
		&
		Binary-like recurrence and scalar-oscillation signatures
		&
		Computational observation
		\\
		
		\([9]^\infty\)
		&
		Ternary-like recurrence and scalar-oscillation signatures
		&
		Computational observation
		\\
		
		\([6]^\infty\)
		&
		Mixed binary-like and ternary-like recurrence and oscillation
		signatures; no single clean recurrence clock
		&
		Compatible with the Chinese-remainder decomposition; detailed
		geometry determined computationally
		\\
		
		\bottomrule
	\end{tabularx}
\end{table}

Thus, constant prime-modulus dynamics is not characterized by stationary
density or entropy. Instead, both scalar observables oscillate as the
system passes through locally recurring templates and globally repeating
\(p\)-adic epochs.

%%%%%%%%%%%%%%%%%%%%%%%%%%%%%%%%%%%%%%%%%%%%%%%%%%%%%%%%%%%

\subsection{A single insertion: effective-seed replacement,
	temporal displacement, and mild densification}

We performed a single non-binary insertion at time \(t=2\), with
\(k\in\{3,4,5,6,7,9\}\), after which all subsequent updates were again
binary. The later evolution retained the characteristic binary
collapse--recovery organization and its successive dyadic epochs.
Thus, a single non-binary update did not establish an autonomous
non-binary regime; instead, it altered the state from which the
subsequent binary evolution proceeded.

\paragraph{Effective-seed replacement} The configuration produced by the non-binary update can therefore be
interpreted as a new effective seed. Subsequent binary epochs are
organized around this transformed template rather than around the
original initial condition. The original seed-return cycle is thereby
replaced by a new binary-type cycle: the large-scale epoch structure
may remain similar, while exact recoveries, when present, concern
descendants of the post-insertion template rather than the original
seed.

\begin{figure}[pos=htbp]
	\centering
	\includegraphics[width=0.7\textwidth]
	%{density_all_k_Moore_neumann_t0_80.png}
	{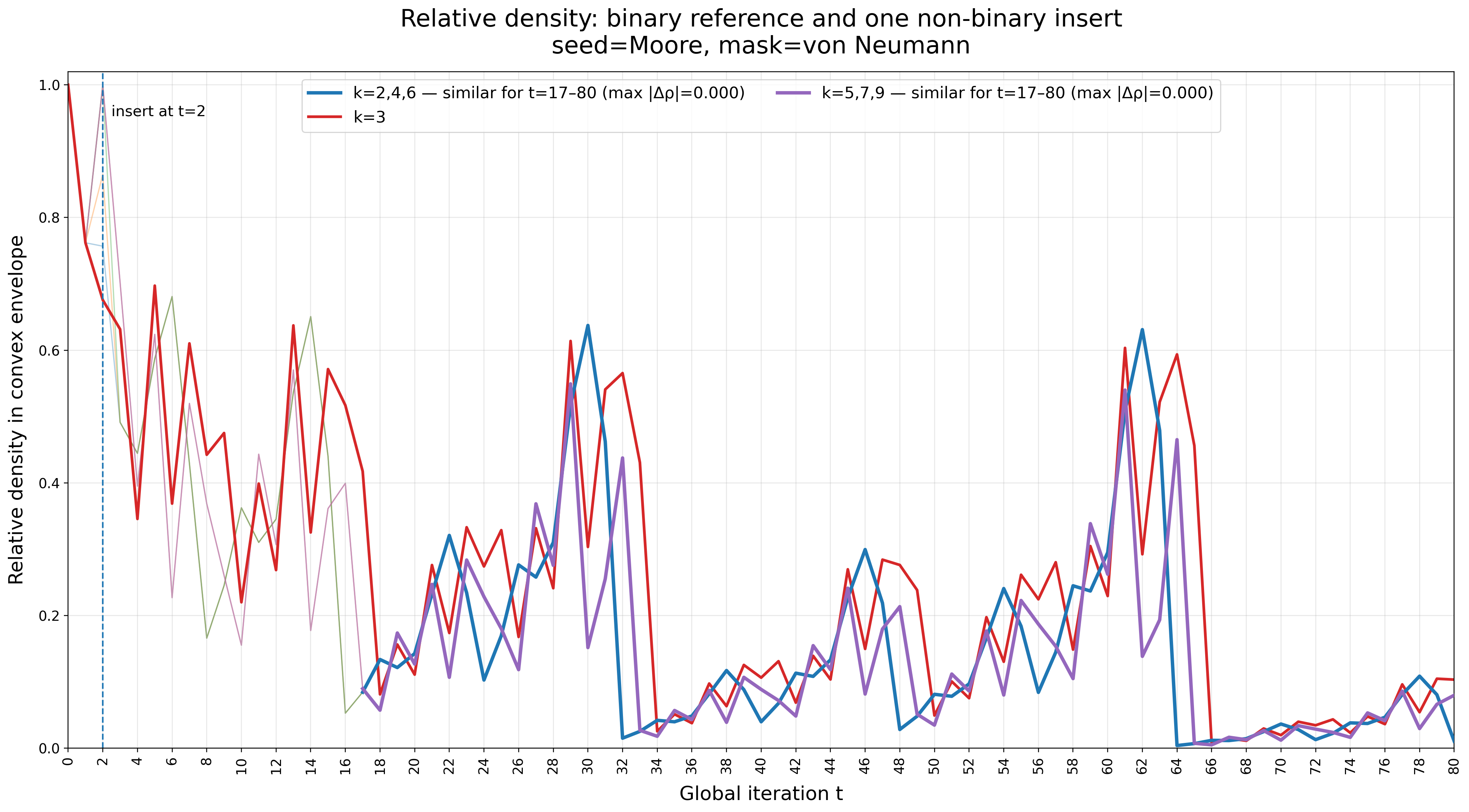}
	\caption{
		Relative-density trajectories after a single non-binary insertion
		at \(t=2\) for the Moore seed and von Neumann mask. After the
		initial transient, the trajectories fall into three phase classes
		of the binary collapse--recovery cycle: \(k=2,4,6\) remain aligned
		with the binary reference, \(k=5,7,9\) are shifted by one
		iteration, and \(k=3\) by two iterations. Within the two
		multi-\(k\) classes, the density trajectories coincide over
		\(t=17,\ldots,80\). Thus, the insertion changes the phase of the
		binary clock without destroying the characteristic binary
		collapse--recovery organization.
	}
	\label{fig:single_insert_phase_classes}
\end{figure}

\paragraph{Binary-clock displacement} This effect depends on
the inserted modulus and the seed--mask pair, the renewed binary
trajectory may remain aligned with the reference or reappear at an
integer phase shift. Figure~\ref{fig:single_insert_phase_classes}
illustrates this particularly clearly for the Moore seed and von
Neumann mask. After the initial transient, the trajectories separate
into three phase classes. The binary reference and the even insertions
\(k=4,6\) remain in the unshifted class; the odd insertions
\(k=5,7,9\) reproduce the same density trajectory shifted by one
iteration, whereas \(k=3\) produces a two-iteration shift. Thus, the
insertion may change the phase of the binary collapse--recovery cycle
without destroying its characteristic temporal organization.

\begin{figure}[pos=htbp]
	\centering
	\includegraphics[width=0.95\textwidth]{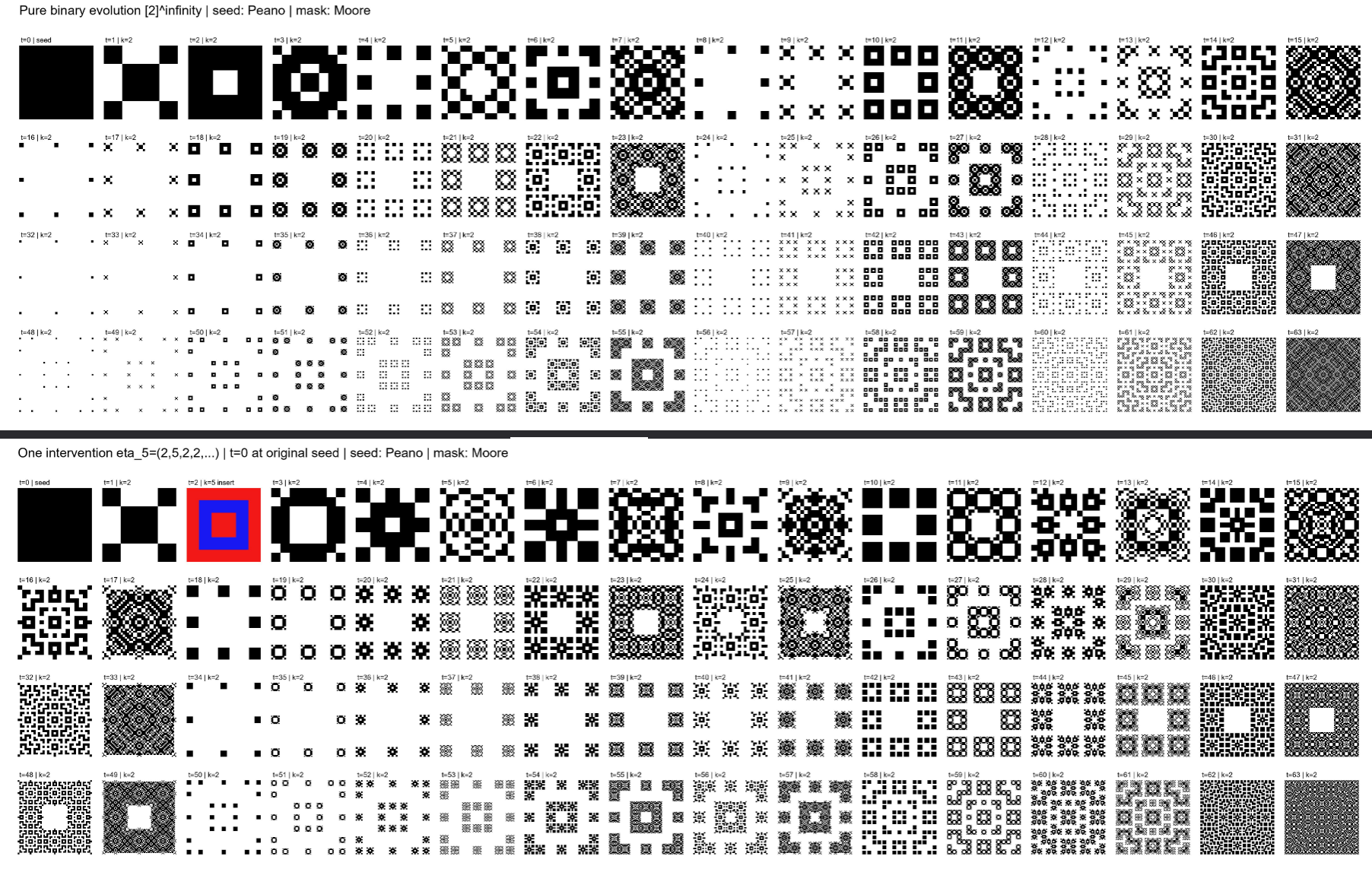}
	\caption{
		A representative case of effective-seed replacement and
		binary-clock displacement after a single non-binary insertion.
		For the Peano seed and Moore mask, the binary evolution
		\([2]^\infty\) (top) is compared with
		\(\eta_{5}=(2,5,2,2,\ldots)\) (bottom). The quinary update at
		\(t=2\) creates a new effective seed from which binary-type
		evolution resumes, with the principal crisis and disappearance
		events shifted by two iterations. The broader dyadic epoch
		organization remains visible, with related template families
		recurring over blocks of lengths \(8,8,16,16,32,\ldots\).
	}
	\label{fig:clock_reset}
\end{figure}

\paragraph{Configuration shift} Figure~\ref{fig:clock_reset} shows the same mechanism directly at the level of configurations. For the Peano seed and Moore mask, a quinary
insertion at \(t=2\) moves the principal collapse and disappearance
events by two iterations relative to the binary reference. The broader
epochal organization nevertheless persists: related template families
continue to recur over blocks of lengths
\(8,8,16,16,32,\ldots\), but they now develop from the
post-insertion configuration.

This separation between large-scale organization and exact state
recovery is consistent with the seed-robust epoch structure shown in
Fig.~\ref{fig:binary_seed_robust_epochs}. Different initial seeds may
populate the same operator-generated scaffold while producing
non-identical complete configurations. A related replacement of the
active template after a non-binary insertion was previously
observed for modular Laplacian dynamics on the hexagonal lattice
\cite{NowakKepczyk2026}.

\paragraph{Mild densification} In addition to the temporal displacement, the odd-modulus insertions produced a modest finite-time increase in mean density over
\(t=2,\ldots,82\). Among these,
the largest mean density gain was obtained for \(k=5\), whereas the
ternary insertion produced the strongest phase displacement.
Temporal density variability remained of comparable magnitude, showing
that the density increase was not accompanied by a transition to a
qualitatively different sustained regime.

A single insertion therefore replaces the effective seed and may
select a different phase of the subsequent binary epoch sequence,
while producing only a modest finite-time density gain. Because all
later updates are binary, the dynamics remains within a renewed
binary-type hierarchy. This motivates the periodic schedules considered
below: can repeated interventions interrupt this renewed binary
organization often enough to sustain an expanding dense regime?

%Fig=10

%Fig=11
%%%%%%%%%%%%%%%%%%%%%%%%%%%%%%%%%%%%%%%%%%%%%%%%%%%%%%%%
%%%%%%%%%%%%%%%%%%%%%%%%%%%%%%%%%%%%%%%%%%%%%%%%%%%%%%%%

\subsection{Periodic non-binary insertions}

A single insertion can replace the effective seed and shift the phase
of the subsequent binary trajectory, but it does not remove the
collapse--recovery organization itself. We therefore ask whether
repeated non-binary interventions can do more: namely, whether they can
prevent repeated binary recovery and sustain a dense expanding regime.

\subsubsection{Short binary tails: from alternating insertions to reduced intervention frequency}
\label{subsec:short-tails}

We first examined the short-tail family
\[
[2,k,2^{\times s}]^\infty,
\]
where \(2^{\times s}\) denotes \(s\) consecutive binary updates.
Thus, \(s=0,1,2\) correspond respectively to
\[
[2,k]^\infty,\qquad
[2,k,2]^\infty,\qquad
[2,k,2,2]^\infty.
\]
These schedules provide three direct cross-sections of the competition
between periodic non-binary intervention and the tendency of the binary
dynamics to recover its sparse dyadic organization.

In many of the odd-modulus schedules considered below, this competition
produced a visually distinct \emph{carpet-like regime} (see Table~\ref{tab:selected_periodic_summary}).

\paragraph{Alternating insertions (\(s=0\)).}
Maximal alternation reveals the clearest parity contrast.
Even moduli retain the characteristic binary scaffold and its
collapse--recovery organization, whereas odd moduli can replace this
pattern by dense carpet-like evolution. The ternary schedule is
particularly distinctive, producing a regular two-phase, sawtooth-like
density oscillation instead of the deep dyadic crises of the binary
reference.

\paragraph{One binary step (\(s=1\)).}
Inserting one binary update between successive non-binary interventions
does not generally restore the binary scaffold. Dense odd-modulus
regimes remain prominent, and the ternary case continues to show
particularly regular persistence. Thus, the first reduction in
intervention frequency is insufficient to recover the sparse binary
organization.

\begin{figure}[pos=htbp]
	\centering
	\includegraphics[width=\textwidth]
	%{density_short_tails_Peano_diag.png}
	{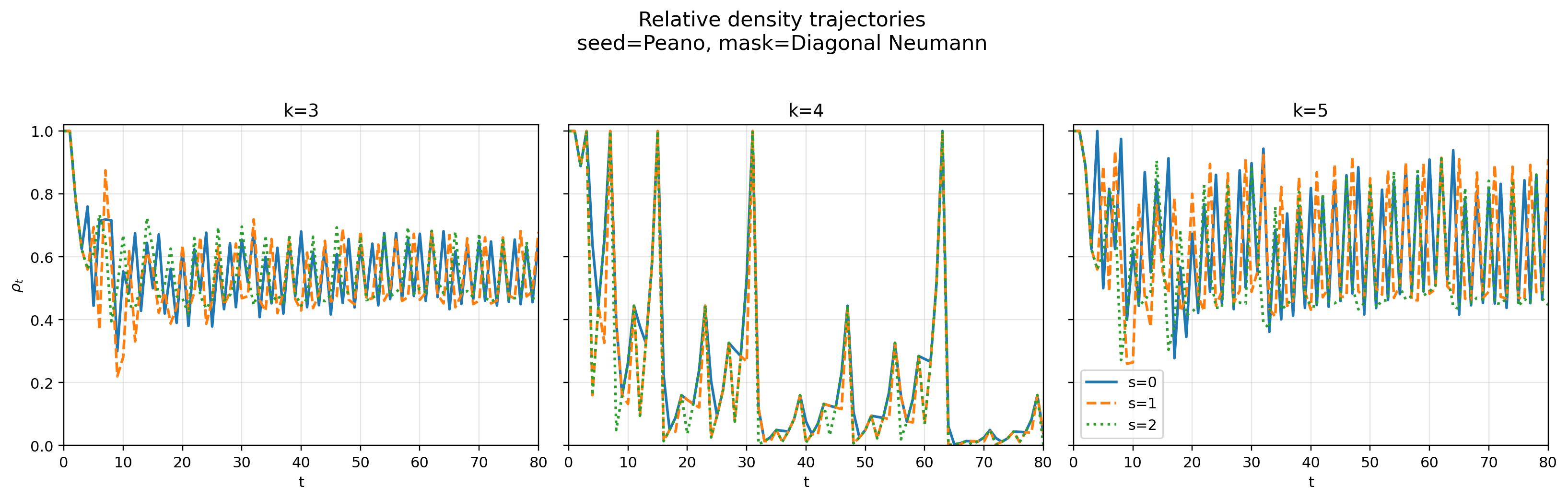}
	\caption{
		\textbf{Relative-density response to increasing binary-tail length.}
		Density trajectories for the Peano seed and diagonal Neumann mask
		under \(k=3,4,5\), with \(s=0,1,2\).
		The ternary schedule \(k=3\) remains in a dense elevated-density
		regime for all three short tails shown. The even control \(k=4\)
		retains the binary-like collapse--recovery scaffold throughout,
		whereas the higher odd modulus \(k=5\) also sustains high density
		but with less regular oscillatory organization than \(k=3\).
		The comparison illustrates how short binary tails sharpen the
		contrast between robust ternary, binary-like even,
		and more geometry-sensitive higher odd insertions.
	}
	\label{fig:short_tail_density}
\end{figure}

\paragraph{Two binary steps (\(s=2\)).}
At \(s=2\), the competition between periodic intervention and binary
recovery becomes more selective rather than disappearing altogether.
In the representative case shown in Fig.~\ref{fig:short_tail_density},
the ternary schedule remains in a dense elevated-density regime,
whereas the even control \(k=4\) remains indistinguishable from the
binary scaffold. The higher odd modulus \(k=5\) also remains dense,
but with more irregular oscillations than \(k=3\). Thus, two intervening
binary steps do not yet universally restore binary organization;
instead, they sharpen the contrast between robust ternary dynamics,
binary-like even dynamics, and less regular higher-odd responses.

The short-tail schedules therefore reveal a progressive competition
between periodic intervention and binary recovery. Alternating
insertions first expose a strong even--odd contrast; one intervening
binary step still permits dense odd-modulus regimes; and by \(s=2\)
the response has become clearly structured by modulus class and
geometry. Across the tested cases, ternary insertion remains the most
robust intervention.

The dependence on binary-tail length is examined systematically in the
following subsection; representative short-tail outcomes are summarized
in Table~\ref{tab:selected_periodic_summary} in the Appendix.

%Fig=12
\begin{figure}[pos=htbp]
	\centering
	%includegraphics[width=0.6\textwidth]{short_tail_cross_sections.png}
	\includegraphics[width=0.6\textwidth]{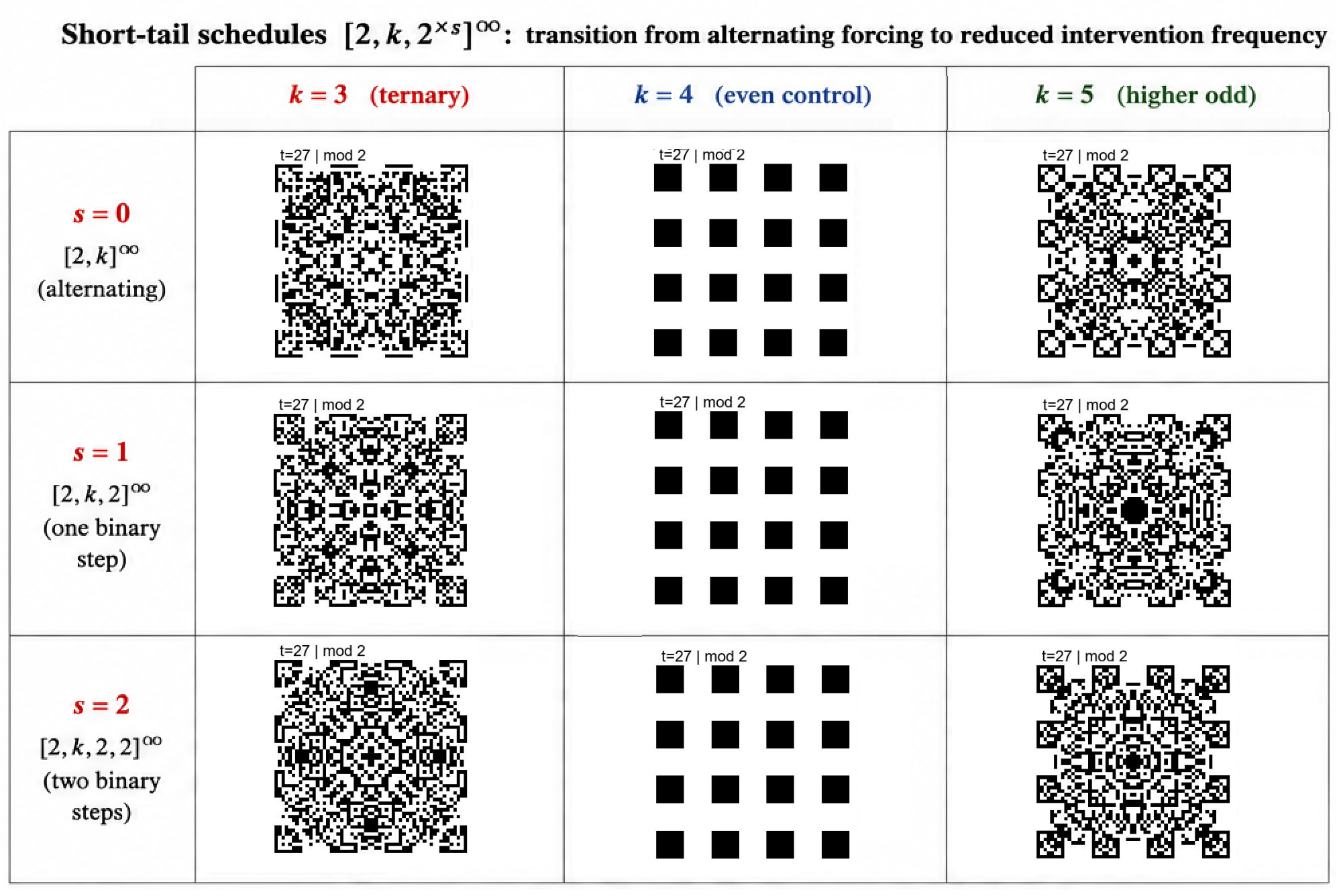}
	\caption{
		\textbf{Short-tail transition under periodic non-binary insertions.}
		Representative configurations for the Peano seed and diagonal
		Neumann mask under three intervention moduli
		(\(k=3,4,5\)) and three binary-tail lengths
		(\(s=0,1,2\)).
		Rows correspond to increasing binary-tail length, while columns
		contrast ternary insertions, an even-modulus control, and a higher odd
		modulus. The ternary schedule remains densely organized across all
		three short tails shown, the even schedule retains a sparse
		binary-like scaffold, and the higher odd modulus occupies an
		intermediate dense but less regular regime.
	}
	\label{fig:short_tail_cross_sections}
\end{figure}

%Fig=13

%Fig=14
\begin{figure}[pos=htbp]
	\centering
	\includegraphics[width=0.8\textwidth]	
	%{density_s2_allk_Peano_diag.png}
	{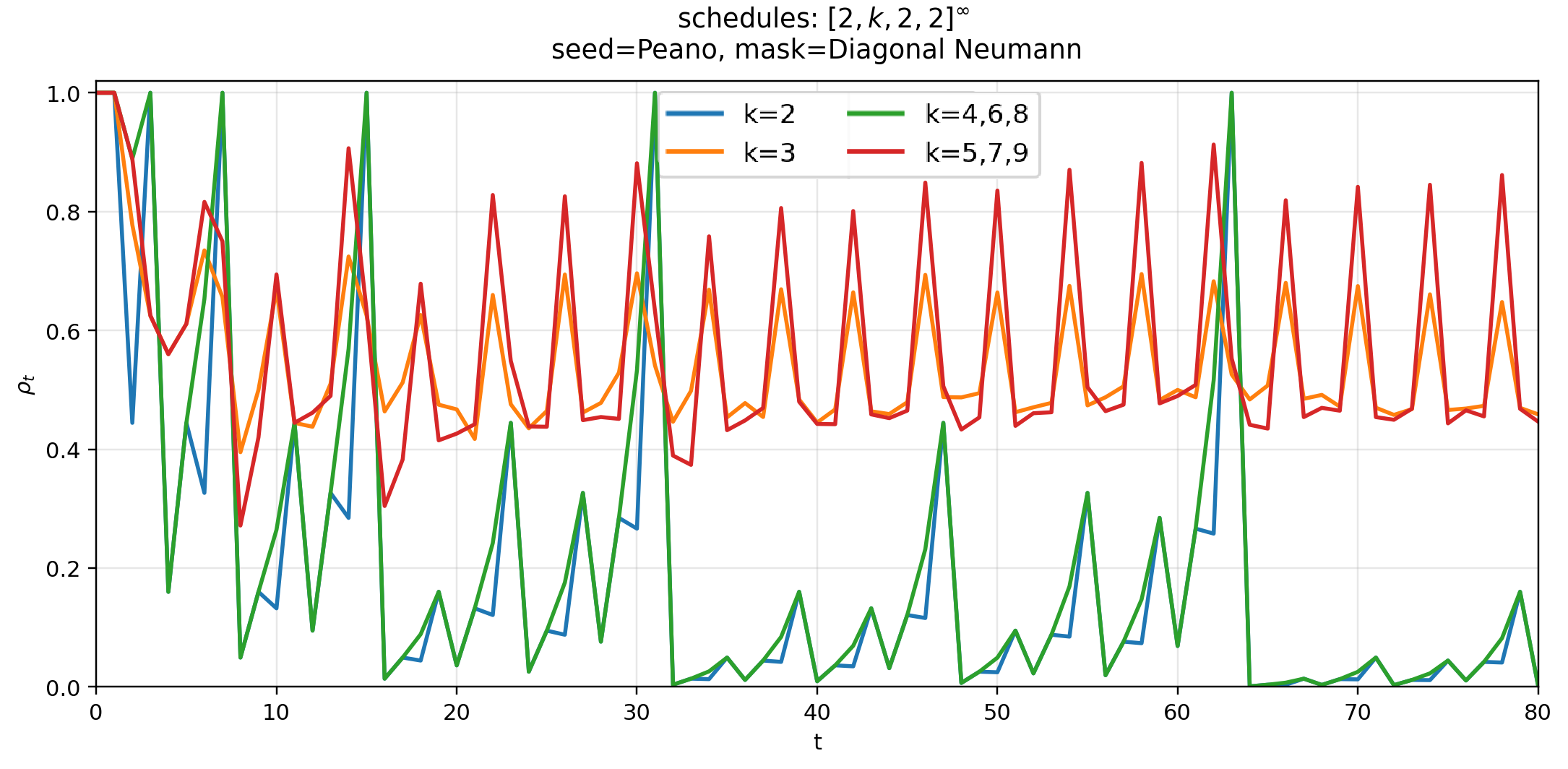}
	\caption{
		\textbf{Modulus classes at \(s=2\).}
		Relative-density trajectories for
		\([2,k,2,2]^\infty\), \(k=2,\ldots,9\), for the Peano seed and
		diagonal Neumann mask. Trajectories differing by at most
		\(10^{-3}\) in relative density over the displayed interval are
		grouped in the legend. The even insertions \(k=4,6,8\) form a
		common binary-like class, while the higher odd moduli
		\(k=5,7,9\) collapse onto a distinct elevated-density class.
		The ternary case \(k=3\) remains separate, emphasizing that the
		short-tail response is structured by modulus class rather than
		varying monotonically with \(k\).
	}
	\label{fig:s2_modulus_classes}
\end{figure}

\paragraph{From shared density phases to exact configuration coincidence.}

The density trajectories shown in Fig.~\ref{fig:s2_modulus_classes} suggested that some higher odd moduli repeatedly enter the same binary phases. This raised a natural question:
does the agreement concern only the global density, or do the systems
actually return to the same spatial configuration? We therefore
compared the complete lattice states at global iterations corresponding
to binary phases shared by the short-tail schedules.

The result was unexpectedly strong. For both the diagonal Neumann and
von Neumann masks, the higher odd moduli \(k=5,7,9\) generated exactly
the same complete configuration at every non-trivial common binary
phase up to \(t=100\), for all three tested seeds and for
\(s=0,1,2\):
\[
u_t^{(5)}=u_t^{(7)}=u_t^{(9)}.
\]
After excluding the initial common binary step at \(t=1\), this
corresponded to exact coincidence at all 33 matched binary checkpoints.
Thus, the previously observed agreement of the density trajectories
was not merely macroscopic: at these phases the systems occupied the
same sites with the same residue values. 

The effect is not universal, however. With the Moore mask, the
higher odd moduli ceased to coincide immediately after the trivial
initial binary phase and retained modulus-specific configurations.
The ternary case \(k=3\) also generally remained outside the
\(k=5,7,9\) equivalence class. A striking exception occurred for the
Peano seed with the von Neumann mask at \(s=2\), where all four odd
moduli \(k=3,5,7,9\) coincided at every matched binary phase.

\paragraph{A parity explanation of the higher-odd coincidence.}
The exact coincidence of the higher odd moduli for the diagonal
Neumann and von Neumann masks has a local explanation. Both
masks have degree \(d_M=4\). If the configuration entering a
non-binary insertion is binary, then at every lattice site
\[
x=(L_Mu)(p)\in\{-4,-3,\ldots,4\}.
\]

\begin{lemma}
	Let \(M\) be a neighborhood mask of degree \(4\), and let
	\(u:\mathbb Z^2\to\{0,1\}\) be a binary configuration. For any odd
	moduli \(k,\ell>4\),
	\[
	[L_Mu]_k\equiv[L_Mu]_\ell\pmod 2
	\]
	pointwise. Consequently, after one subsequent binary update,
	\[
	\bigl[L_M[L_Mu]_k\bigr]_2
	=
	\bigl[L_M[L_Mu]_\ell\bigr]_2 .
	\]
\end{lemma}

\noindent
Indeed, for \(x\in[-4,4]\) and odd \(k>4\),
\[
[x]_k=
\begin{cases}
	x, & x\geq0,\\
	k+x, & x<0.
\end{cases}
\]
Hence the parity of \([x]_k\) is independent of the particular odd
modulus \(k>4\):
\[
[x]_k\bmod 2=
\begin{cases}
	x\bmod2, & x\geq0,\\
	(x+1)\bmod2, & x<0.
\end{cases}
\]
Thus the configurations produced by different higher odd moduli may
contain different residue values immediately after the insertion, but
they have the same pointwise parity. Since the Laplacian is an
integer-linear operator, the next reduction modulo \(2\) maps all of
them to exactly the same binary configuration.
\[
\begin{array}{c|rrrrrrrrr}
	x & -4 & -3 & -2 & -1 & 0 & 1 & 2 & 3 & 4 \\ \hline
	\left[x\right]_5 \bmod 2 & 1 & 0 & 1 & 0 & 0 & 1 & 0 & 1 & 0 \\
	\left[x\right]_7 \bmod 2 & 1 & 0 & 1 & 0 & 0 & 1 & 0 & 1 & 0 \\
	\left[x\right]_9 \bmod 2 & 1 & 0 & 1 & 0 & 0 & 1 & 0 & 1 & 0
\end{array}
\]

More generally, the
same mechanism predicts a common higher-odd synchronization class for
odd \(k>d_M\). This also explains why the \(5,7,9\) coincidence need
not persist for the Moore mask, whose degree is \(8\), and suggests
that \(k=9,11,13,\ldots\) should form the corresponding class there.
Below this threshold, additional coincidences may still occur for
particular configurations, as suggested by the exceptional involvement
of \(k=3\) observed in one seed--mask--schedule combination.

%Fig=15
\begin{figure}[pos=htbp]
	\centering
	\includegraphics[width=0.75\textwidth]{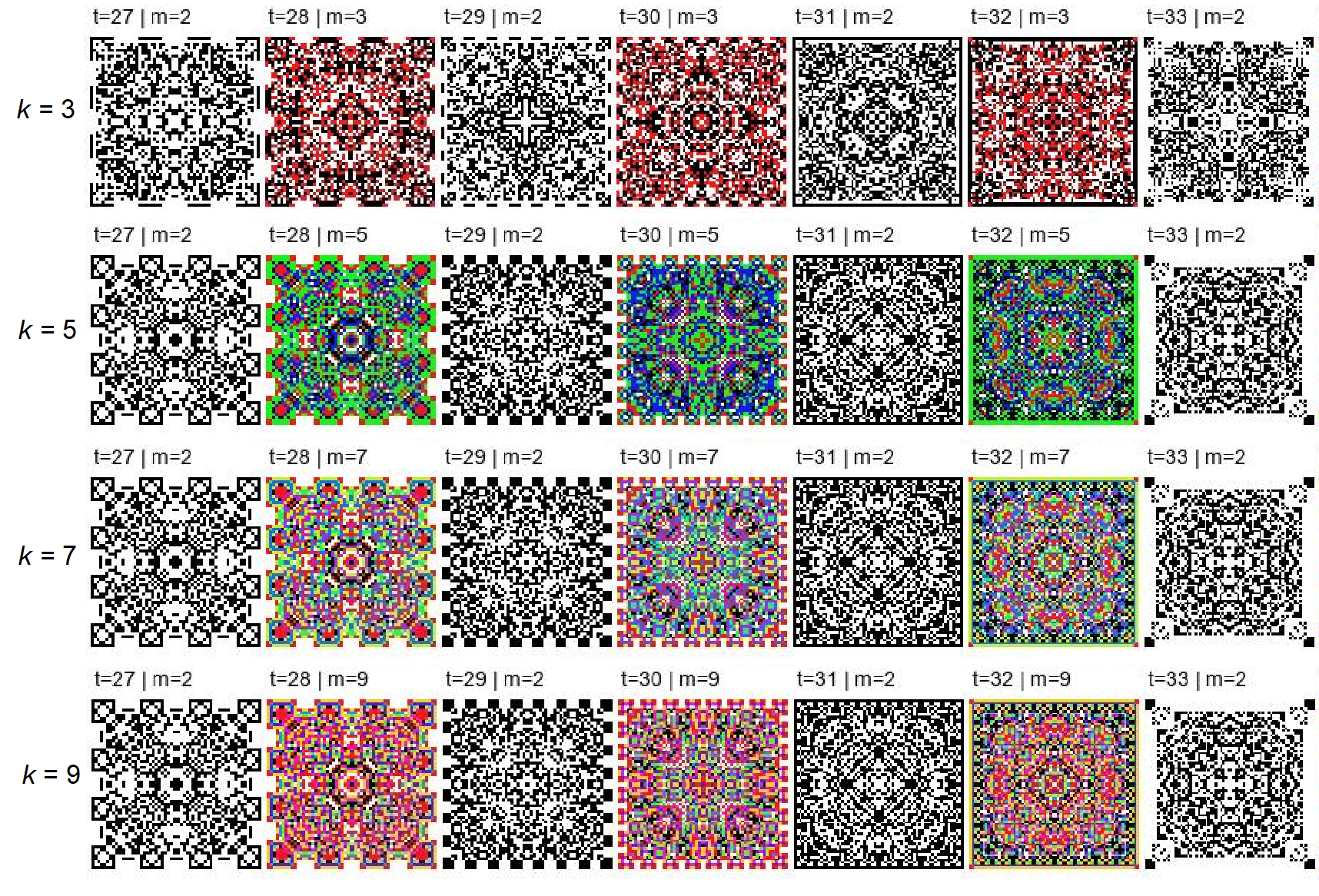}
\caption{
	\textbf{From shared density phases to identical patterns.}
	For the Peano seed and diagonal Neumann mask, the common density phases
	of \(k=5,7,9\) correspond to exact coincidence of the complete lattice
	state, while the ternary case remains distinct.
}
	\label{fig:the_same}
\end{figure}

%%%%%%%%%%%%%%%%%%%%%%%%%%%%%%%%%%%%%%%%%%%%%%%%%%%%%%%%%%
%%%%%%%%%%%%%%%%%%%%%%%%%%%%%%%%%%%%%%%%%%%%%%%%%%%%%%%%%%
%%%%%%%%%%%%%%%%%%%%%%%%%%%%%%%%%%%%%%%%%%%%%%%%%%%%%%%%%%
%%%%%%%%%%%%%%%%%%%%%%%%%%%%%%%%%%%%%%%%%%%%%%%%%%%%%%%%%%

\subsubsection{Tail-length sensitivity in
\([2,k,2^{s}]^\infty\)}

\label{subsec:tail-sensitivity}

We next extended the binary tail in
\[
[2,k,2^{s}]^\infty,
\]
where \(2^{s}\) denotes \(s\) consecutive binary updates.
Successive non-binary interventions are therefore separated by
\(s+1\) binary updates. Density was evaluated in the main window
\[
W=\{33,\ldots,200\}.
\]

\begin{figure}[pos=htbp]
	\centering
	\includegraphics[width=\textwidth]
	%{Fig_B_small_seed_median_heatmap.png}
	{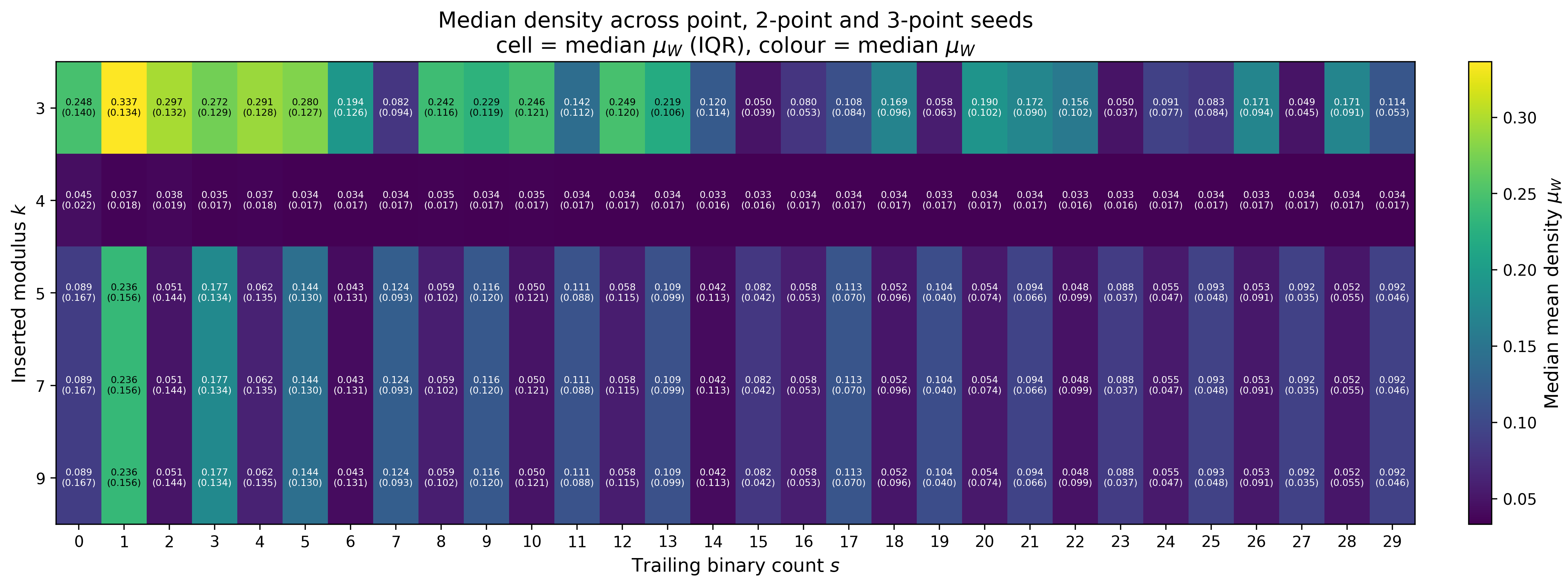}
	\caption{
		\textbf{Tail-length sensitivity across inserted moduli.}
		Median mean density in \(W=[33,200]\) across three small seeds for the
		diagonal Neumann mask. Ternary insertions maintain the broadest
		high-density region, the even case remains sparse, and the higher odd
		moduli show strongly oscillatory responses to binary-tail length.
	}
	\label{fig:tail_heatmap}
\end{figure}

To expose schedule effects with minimal geometric overlap, the
systematic scan was performed for three small seeds under the diagonal
Neumann mask. Figure~\ref{fig:tail_heatmap} shows the median mean
density across these seeds. Increasing the binary tail did not produce
a monotone loss of density. Instead, favourable and unfavourable tail
lengths alternated, revealing a strong phase dependence of the
periodic insertions.

The inserted moduli separated into distinct response classes.
Ternary insertions formed the broadest high-density branch, whereas the
even case \(k=4\) remained close to the sparse binary regime. The
higher odd moduli \(k=5,7,9\) showed lower and more strongly
oscillatory profiles. Their aggregate responses were undistinguishable, following the lemma.

%Fig=16

For \(k=3\), the non-monotone response showed pronounced minima near
\[
s=7,\;15,\;23,\;31,
\]
corresponding to \(s+1=8,16,24,32\) binary updates between successive
non-binary insertions. Their regular spacing links the density
losses to the characteristic binary epoch structure. Tail length
therefore acts not simply as a measure of intervention frequency, but
as a phase-selecting control parameter.

The tail-length scan therefore extends the short-tail result in an
important way. Dense regimes are not controlled solely by how often a
non-binary update is applied. Their persistence also depends on
\emph{when} that update meets the evolving binary hierarchy. Ternary
insertion is the most robust tested intervention, but even for \(k=3\)
the binary clock contains recurrent phases at which densification is
strongly suppressed.

%%%%%%%%%%%%%%%%%%%%%%%%%%%%%%%%%%%%%%%%%%%%%%%%%%%%%%%%%%
%%%%%%%%%%%%%%%%%%%%%%%%%%%%%%%%%%%%%%%%%%%%%%%%%%%%%%%%%%

\subsubsection{Seed extent attenuates tail-length sensitivity}
\label{subsec:seed-size}

The pronounced oscillations above were obtained for small seeds with
the diagonal Neumann mask, a setting in which schedule-dependent
differences were particularly clear. We therefore asked whether the
same phase-sensitive response persists when the initial configuration
is more extended.

\begin{figure}[pos=htbp]
	\centering
	\includegraphics[width=0.99\textwidth]
	%{Fig_seed_size_heatmap_mean_density.png}
	{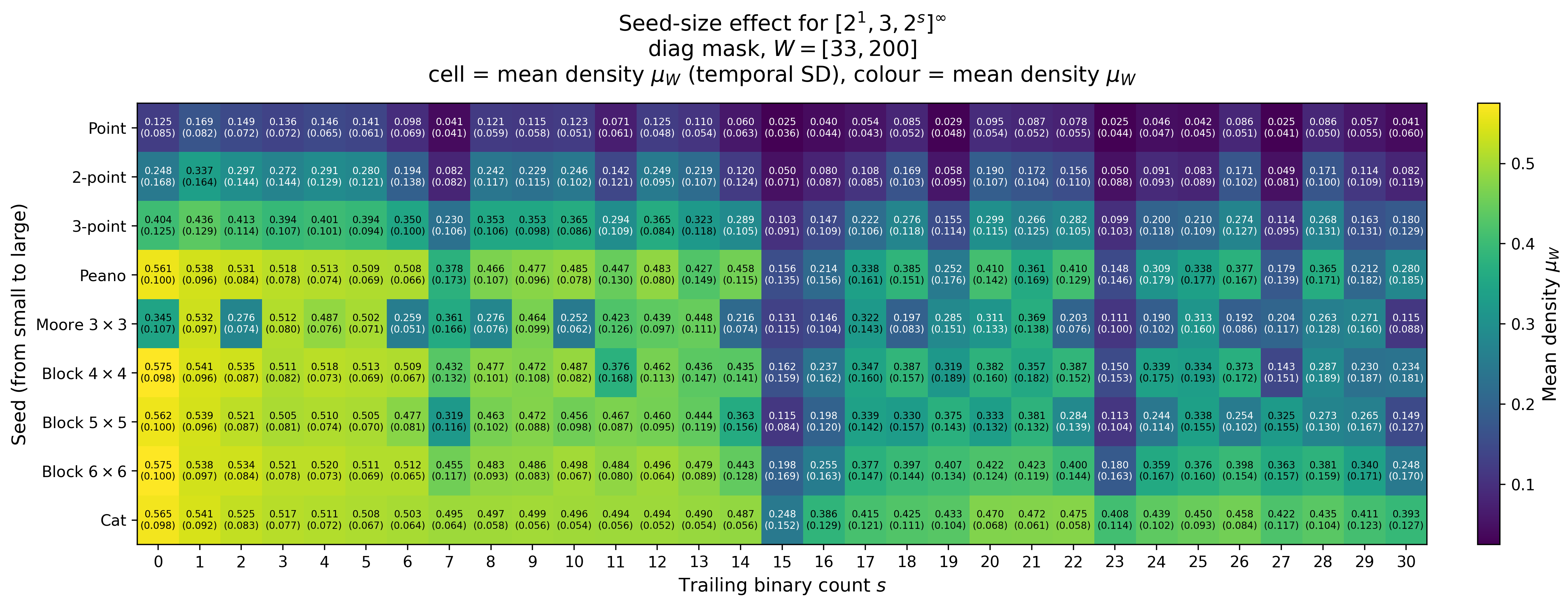}
	
	\vspace{4mm}
	
	\includegraphics[width=0.99\textwidth]{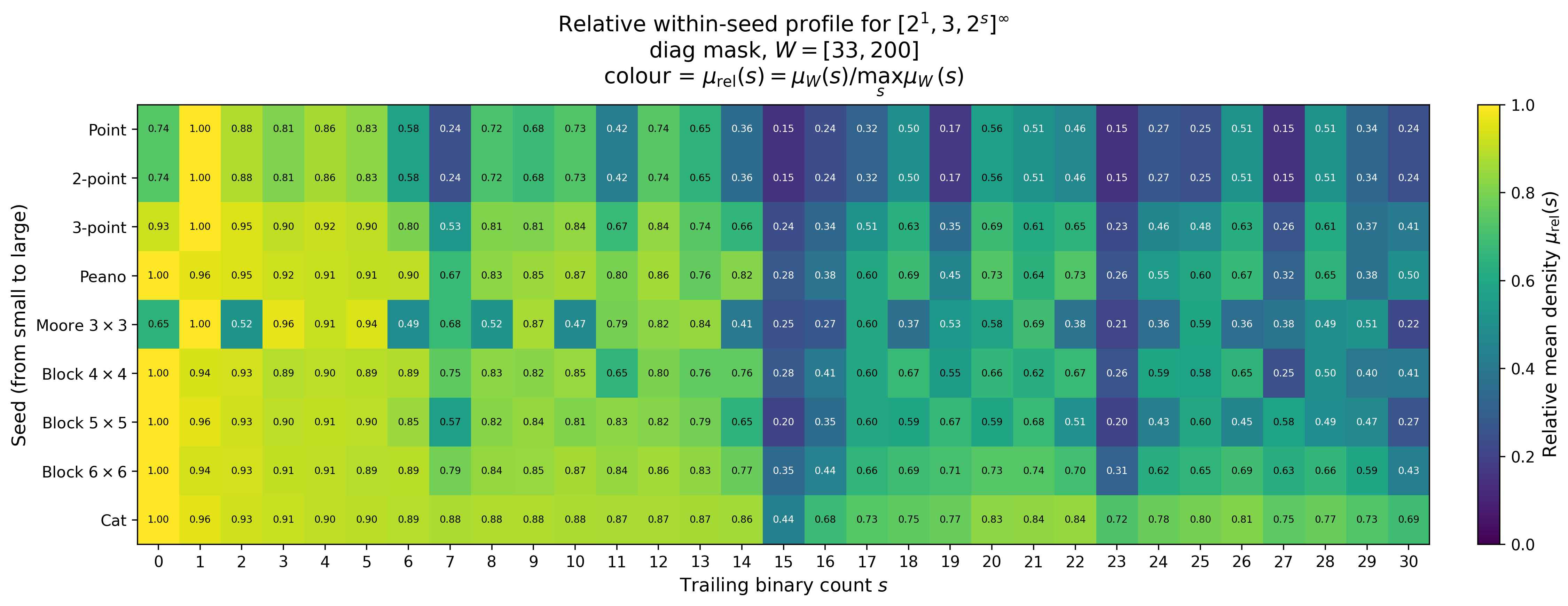}
	%{Fig_seed_size_heatmap_relative_density.png}

\caption{\textbf{Attenuation of tail-length sensitivity for extended seeds.}
	Mean density for \(k=3\) and the diagonal Neumann mask in
	\(W=[33,200]\). The upper panel shows absolute density
	\(\mu_W(s)\); the lower panel shows each seed relative to its own
	maximum,
	\[
	\mu_{\mathrm{rel}}(s)=
	\frac{\mu_W(s)}{\max_s \mu_W(s)}.
	\]
	Small seeds exhibit pronounced minima at
	\(s=7,15,23,31\), whereas extended seeds generally retain higher
	relative density at these unfavourable tail lengths.}
	\label{fig:seed_size_heatmaps}
\end{figure}

Figure~\ref{fig:seed_size_heatmaps} shows the result for \(k=3\).
Extended seeds generally produced higher absolute density, but the
more important effect was a reduction in the amplitude of the
tail-length response. Small seeds suffered deep density losses at
unfavourable values of \(s\), whereas extended seeds retained a larger
fraction of their own maximum density across the same tail lengths.

Thus, seed extent acts as a buffer against unfavourable perturbation
phases. A larger initial support generates overlapping descendants for
longer, so that binary separation is less able to expose the deep
density troughs seen for small seeds. The attenuation is not strictly
ordered by support size, however: Peano, Moore, compact blocks, and
image seeds retain distinct profiles, showing that seed geometry
matters in addition to extent.

%Fig=17

Seed geometry therefore modulates the phase sensitivity identified
above. Small seeds provide the clearest expression of the binary
tail-length structure, whereas more extended seeds progressively
attenuate the same response. Tail length determines when the
intervention encounters the binary hierarchy, while seed geometry
determines how strongly that phase is expressed in the resulting
density.

%%%%%%%%%%%%%%%%%%%%%%%%%%%%%%%%%%%%%%%%%%%%%%%%%%%%%%%%%%
%%%%%%%%%%%%%%%%%%%%%%%%%%%%%%%%%%%%%%%%%%%%%%%%%%%%%%%%%%

\section{Discussion and conclusions}

The world of constant-modulus Laplacian trajectories is surprisingly
orderly. Under a fixed prime modulus, the evolution is tied to a strict
\(p\)-adic clock. In the binary case, a seed grows through successive
paired shape epochs, whose characteristic durations double, until the
configuration collapses into spatially separated seed copies and a new,
larger epoch begins. The ternary system tells the same story in threes:
three characteristic developmental phases are followed by renewed
separation and rebuilding on the next \(3\)-adic scale. Other prime
moduli follow the same principle, each with its own \(p\)-adic clock.
The seed and neighborhood mask may change the geometry, the amount of
overlap, and the visible detail, but the underlying temporal scaffold
is remarkably difficult to escape.

Against this highly organized background, the effect of changing the
modulus in time becomes especially striking.

The results distinguish two qualitatively different effects of
non-binary insertions into binary modular Laplacian dynamics. A single
insertion changes the active configuration and may shift the subsequent
binary replication clock, but it does not destroy the underlying
collapse--recovery organization. In this sense, an isolated insertion
acts mainly through effective-seed replacement: after the intervention,
the binary operator continues to organize the dynamics, although from
a modified state and possibly from a different phase of its dyadic
epoch structure.

Periodic insertions produce a substantially different response.
For suitable odd moduli and binary-tail lengths, repeated recovery of
the sparse binary scaffold is suppressed and replaced, over the
investigated horizon, by long-lived densely occupied carpet-like
regimes. The response is not ordered simply by the magnitude of the
inserted modulus. Even insertions remain predominantly binary-like,
whereas odd insertions can sustain dense states, with \(k=3\) producing
the broadest and most robust high-density response among the tested
moduli. Thus, the effect of a non-binary update depends not only on its
modulus, but also on the state of the binary dynamics at which it is
applied.

This phase dependence is particularly evident in the binary-tail
experiment. Density varies non-monotonically with \(s\), and for
\(k=3\) pronounced minima occur near
\[s=7,15,23,31,\]
corresponding to \(8,16,24,32\) binary updates between successive
insertions. Their regular spacing suggests an interaction with the
underlying binary epoch structure. Tail length therefore acts not
merely as a measure of insertion frequency, but as a phase-selecting
parameter. The same non-binary update may be dynamically effective at
one phase and nearly ineffective at another. This observation suggests
that the relevant object is not the inserted modulus alone, but the
pair consisting of the current configuration and the position of the
intervention within the binary replication hierarchy.

One of the most unexpected results is the exact convergence of
trajectories generated by different higher odd moduli. For the
diagonal Neumann and von Neumann masks, the systems with
\(k=5,7,9\) differ during the non-binary step but repeatedly return
to exactly the same complete configuration at the following binary
phases.
The lemma proved here explains this coincidence. 

Seed geometry provides another level of control. Extended seeds
generally attenuate the response to unfavourable tail lengths, whereas
small seeds expose the binary phase structure most clearly. The
observable dynamics therefore reflects an interaction between the
temporal schedule, the neighborhood operator, and the geometry of the
evolving state.

The periodic family
\[
[2,k,2^{\times s}]^\infty
\]
is only a small subset of possible time-dependent modular schedules.
More general sequences
\[
\mathbf{k}=(k_1,k_2,k_3,\ldots)
\]
raise a broader classification problem: which properties of the
schedule determine recurrence, densification, collapse, or convergence
of distinct trajectories? The same question may be extended to
non-binary initial configurations, where support geometry and residue
information become independent sources of dynamical memory.

This leads naturally to an inverse problem: given a prescribed
geometric or dynamical target, determine the initial configuration,
neighborhood mask, and modular schedule that reproduce it, exactly or
approximately. The present experiments identify several relevant
control variables---modulus class, binary-tail phase, mask degree, and
seed geometry---but explore only a restricted region of this larger
parameter space.

This inverse viewpoint connects the present problem with broader
questions in nonlinear dynamics and pattern formation, where controlled
perturbations are used to select, stabilize, or redirect emerging
spatiotemporal structures \cite{Burkart2022,Gershenson2025}.
Here the control variable is discrete and symbolic: the modular
schedule acts as a sequence of interventions applied to an evolving
spatial state. Modular Laplacian dynamics may therefore provide a
simple setting for studying inverse design and control of self-organized
patterns, related in spirit to recent work on programmable
spatiotemporal pattern formation in cellular-automaton-like systems
\cite{Richardson2024}.

These possibilities remain open beyond the finite computational
horizon studied here, and no asymptotic persistence or aperiodicity is
claimed. The parity mechanism already explains the higher-odd
synchronization above the mask-degree threshold, while the origin of
the regularly spaced tail-length minima and the
configuration-dependent coincidences below that threshold remains open.

\section*{Data availability}

Representative datasets and supplementary materials generated
during the study are available via Figshare:\\
DOI: 10.6084/m9.figshare.32592201

\paragraph{Supplementary material.}
The supplementary material contains the Python source code used in the
computational experiments together with the input seed configurations.

\newpage

\subsection*{Acknowledgments}

The author thanks the anonymous reviewers for valuable comments and suggestions.

\subsection*{Disclosure statement}

No potential conflict of interest was reported by the author.

\section*{Declaration of generative AI and AI-assisted technologies
	in the manuscript preparation process}

During the preparation of this work, the author used ChatGPT (OpenAI)
to assist with language editing, manuscript restructuring, and the
review and debugging of selected Python code. The author critically
reviewed, verified, and edited all AI-assisted output as needed and
takes full responsibility for the content of the publication.

\printcredits

\begin{table}[H]
	\centering
\caption{\textbf{Appendix table A}. Phase shifts of post-insertion density
	trajectories relative to the binary reference, estimated over
	\(t=17--80\). An entry \(q(+r)\) denotes a shift of \(r\) iterations
	after a modulo-\(q\) insertion; 0 denotes no shift, and braces group
	moduli with the same shift. See Figs.~\ref{fig:clock_reset} and
	\ref{fig:single_insert_phase_classes}.}
	\label{tab:single_insert_lags}
	
	\scriptsize
	\setlength{\tabcolsep}{5pt}
	\renewcommand{\arraystretch}{1.25}
	
	\begin{tabularx}{\textwidth}{
			>{\raggedright\arraybackslash}p{0.13\textwidth}
			>{\raggedright\arraybackslash}X
			>{\raggedright\arraybackslash}X
			>{\raggedright\arraybackslash}X
		}
		\toprule
		\textbf{Seed}
		& \textbf{Diagonal Neumann mask}
		& \textbf{Moore mask}
		& \textbf{von Neumann mask} \\
		\midrule
		
		Point
		& \(3(+2)\); \(\{5,7,9\}(0)\)
		& \(\{3,5,7\}(+2)\); \(9(0)\)
		& \(3(+2)\); \(\{5,7,9\}(0)\) \\
		
		Peano
		& \(\{3,5,7,9\}(+2)\)
		& \(\{3,5\}(+2)\); \(\{7,9\}(+1)\)
		& \(\{3,5,7,9\}(+1)\) \\
		
		Moore
		& \(3(+2)\); \(\{5,7,9\}(+1)\)
		& \(3(+1)\); \(\{5,7,9\}(+2)\)
		& \(3(+2)\); \(\{5,7,9\}(+1)\) \\
		
		\bottomrule
	\end{tabularx}	
\end{table}

\begin{table}[H]
	\centering
\caption{\textbf{Appendix table B}. Density and resistance to collapse
	for selected \(s=2\) schedules in \(W=[33,200]\), for the Peano seed
	and diagonal Neumann mask. Relative density is measured within the
	convex bounding polygon. \(Q_{0.10}(\rho)\) is the 10th percentile of
	the density trajectory; grouped classes are summarized by medians
	across moduli. See Fig.~\ref{fig:short_tail_density}.}
	\label{tab:selected_periodic_summary}
	
	\begin{tabular}{lcccc}
		\toprule
		Class & Schedule & \(\mu_W\) & Gain vs. binary & \(Q_{0.10}(\rho)\) \\
		\midrule
		Binary reference & $[2]^\infty$ & 0.066 & reference & 0.002 \\
		Ternary & $[2,3,2,2]^\infty$ & 0.531 & +710\% & 0.471 \\
		Even class & $k=4,6,8$ & 0.076 & +16\% & 0.002 \\
		Higher odd class & $k=5,7,9$ & 0.575 & +777\% & 0.456 \\
		\bottomrule
	\end{tabular}

\end{table}

\bibliographystyle{cas-model2-names}

\bibliography{references}

\end{document}